\pdfoutput=1
\documentclass[a4paper]{article}
\usepackage[latin9]{inputenc}
\usepackage[english]{babel}
\usepackage{amstext,amsmath,amsfonts,amssymb,amsbsy,bm,amsthm}
\usepackage{graphicx}
\usepackage[usenames,dvipsnames]{color}
\usepackage[]{subfig}
\usepackage[normalem]{ulem}
\usepackage[bottom=3cm,top=2cm,left=2cm,right=2cm]{geometry}

\numberwithin{equation}{section}

\newtheorem{definition}{Definition}[section]
\newtheorem{remark}[definition]{Remark}
\newtheorem{prob}{Problem}

\title{Optimal-order isogeometric collocation at Galerkin superconvergent points}

\author{M. Montardini\footnote{Dipartimento di Matematica, Universit\`a degli Studi di Pavia, Italy \texttt{monica.montardini01@universitadipavia.it}}, 
  	    G. Sangalli\footnote{Dipartimento di Matematica, Universit\`a degli Studi di Pavia, Italy, e
          Istituto di Matematica Applicata e Tecnologie Informatiche ``E. Magenes'' del CNR, Pavia, Italy
          \texttt{giancarlo.sangalli@unipv.it}}, 
 		L. Tamellini\footnote{Istituto di Matematica Applicata e Tecnologie Informatiche ``E. Magenes'' del CNR, Pavia, Italy
          \texttt{tamellini@imati.cnr.it}}}

\begin{document}

\maketitle

\paragraph{Abstract}
In this paper we investigate numerically the order of convergence of 
an isogeometric collocation method that builds upon
the least-squares collocation method presented in \cite{anitescu2015isogeometric} 
and the variational collocation method presented in 
\cite{gomez2016variational}. The focus is on smoothest B-splines/NURBS
approximations, i.e, having global $C^{p-1}$ continuity for polynomial degree $p$.  
Within the framework of \cite{gomez2016variational}, we select  as collocation points
a subset of those considered in
\cite{anitescu2015isogeometric}, which are related to the Galerkin
superconvergence theory. With our choice,   that features local
symmetry of the collocation stencil, we improve the 
convergence behaviour with respect to  \cite{gomez2016variational},
achieving   optimal $L^2$-convergence for odd degree 
B-splines/NURBS approximations. 
The same optimal order of convergence  is seen in
\cite{anitescu2015isogeometric}, where, however a  least-squares formulation 
is adopted. 
Further careful study is needed, since the
robustness of the method and its mathematical foundation are still unclear.


\section{Introduction}

The splines-based collocation method for solving differential equations has about fifty years of history. 
The first references are \cite{bickley1968piecewise,fyfe1969use}, where cubic $C^2$ splines are
used to solve a second order two-point boundary value problem. In particular, in
order to achieve optimal convergence, \cite{fyfe1969use} collocates  a modified 
equation, where the modification is obtained by constructing a  suitable interpolant
of the true solution. An extension of this approach to multivariate
(tensor-product) splines and partial differential equations is
studied in \cite{houstis1988convergence}, while extensions to $m$-order differential
equations are found in \cite{russell1972collocation} and in particular
in \cite{de1973collocation}, where the optimality of the
method is achieved  by collocating the original, unperturbed, equation at
suitably selected collocation points, i.e, Gaussian quadrature
points. 
The method {only} works for splines of continuity $C^{m-1}$ and
degree $m+k-1$, with $k\geq m$. Splines-based collocation has been 
successfully applied also to integro-differential equations on curves, and
to  the boundary element method for planar domains  (see
\cite{arnold1985convergence} and references therein). 

The interest and development of splines-based collocation methods 
for partial differential equations has been
driven in the last decade by isogeometric analysis (see \cite{Hughes2005,IGA-book,Auricchio:2010,
Auricchio2013113,Collocation2,da2012avoiding,kiendl2015isogeometric,reali2015isogeometric,
de2015isogeometric,reali2015introduction,casquero2016isogeometric,matzen2013point,
gomez2014accurate,manni2015isogeometric,gomez2016variational}
and references therein). The motivation is computational
efficiency:  isogeometric collocation is more efficient than the isogeometric 
Galerkin method, at least for standard code implementations, 
see \cite{Schillinger:2013}. In particular, the assembly of
system matrices is much faster for collocation than
for Galerkin (unless one adopts recent quadrature algorithms as in
\cite{calabro2016fast}). On the other hand, contrary to the Galerkin method,
isogeometric collocation based on maximal regularity splines 
has always been reported suboptimal in literature, when the error is measured in 
$L^2$ or $L^\infty$ norm. For example, the $L^2$ norm of the error
of the collocation method at Greville points, studied in
\cite{Auricchio:2010}  for a second-order  elliptic problem, 
converges under $h$-refinement as $O(h^{p-1})$ or
$O(h^{p})$, when the degree $p$ is odd or even, respectively,
while the optimal interpolation error is $O(h^{p+1})$ 
regardless of the parity of $p$ for a smooth solution. We remark that the previous ideas of
\cite{fyfe1969use,de1973collocation} cannot be applied directly  to
the isogeometric case since \cite{fyfe1969use} would require a complex modification
of the equation (this approach however deserves further investigation)
and \cite{de1973collocation} does not work for maximal smoothness
splines, which represent the most interesting choice in this framework.

Collocating the equation at Greville points (obtaining the method to which we refer here as 
Collocation at Greville Points, C-GP) is a
common choice since Greville points are classical interpolation points for arbitrary 
degree and regularity splines, well studied in literature, see e.g. \cite{DeBoor}. 
There is however an alternative and  interesting approach, 
from  \cite{anitescu2015isogeometric} and
\cite{gomez2016variational}.
In particular, \cite{gomez2016variational}  introduces an ideal collocation
scheme whose solution coincides with the solution of the  Galerkin
method, thus recovering optimal convergence. This scheme uses as
collocation points the so-called Cauchy-Galerkin points, a well-chosen
subset of  the zeros of the Galerkin residual.  These
points are not known a-priori, and {therefore} 
\cite{gomez2016variational} selects as approximated
Cauchy-Galerkin  points  
the points where, under some hypotheses (we will return
on this point later on, in Section~\ref{sec:sppoints}), one can prove
superconvergence of the second derivatives of the Galerkin
solution.   Indeed, for a Poisson problem the residual is equivalent to the error on the
approximation of the second derivatives. This is an idea from the
previous paper \cite{anitescu2015isogeometric}: 
if we constrain the numerical residual to be zero where the Galerkin residual is estimated
to be zero up to higher order terms, then the computed numerical solution is
expected to be close to the Galerkin numerical solution up 
to higher order terms as well. 
There are however two difficulties: the first and most
relevant one is that also the superconvergent points are not known with
enough accuracy everywhere in the computational domain;
the second one is that there are more Galerkin
superconvergent points than degrees-of-freedom,  $n_{dof}$, for
maximal smoothness splines (the superconvergent points are about
$2n_{dof}$). Indeed, \cite{anitescu2015isogeometric} proposes to compute
a solution of the overdetermined linear system by a least-square
approximation. This approach, which is more expensive  than
collocation, achieves optimal convergence for odd degrees and
one-order suboptimal for even degrees. We refer to it  as Least-Squares approximation  at
Superconvergent Points (LS-SP). Instead, \cite{gomez2016variational} 
designs a well-posed collocation scheme by selecting only $n_{dof}$ collocation points
among those used in \cite{anitescu2015isogeometric}. Roughly speaking,
one superconvergent point per element is used as collocation point,
i.e., one every other superconvergent point (as shall
be clearer later), and therefore in this paper we denote this method as Collocation at
Alternating Superconvergent Points (in short C-ASP). The $L^2$
convergence of C-ASP is one-order suboptimal for any degree,
i.e., the $L^2$-error decays as $O(h^{p})$ for any $p$, which means that
the lack of accuracy in the estimated location of the  superconvergent
points  affects the  convergence  behaviour of the collocation method
C-ASP. 

However, we have an interesting and useful finding to report in this
paper. In the framework of \cite{gomez2016variational},  we
propose a new criterion for selecting the subset of 
superconvergent points, which features 
local symmetry and gives improved convergence properties {compared to} C-ASP. 
Roughly speaking, we {propose to} take two (symmetric) superconvergent points 
in every other element.  This method, which we refer to as  Clustered Superconvergent Points (C-CSP), 
features the same convergence order as the LS-SP approach,
i.e., optimal convergence for odd degrees {in $L^2$ and $L^\infty$ norm}.  
Thus, we finally have an optimally convergent isogeometric collocation scheme
with cubic $C^2$ splines (see \cite{Schillinger:2013} for a discussion
on the relevance of this case).

The results we have obtained are preliminary and, while
some ``magical'' error cancellation  happens with the C-CSP collocation
point selection,  perhaps due to the local symmetry of the collocation
stencil, we are still unable to provide a rigorous
convergence proof for C-CSP (nor for LS-SP or C-ASP).
Furthermore, we have considered quite simple
numerical benchmarks,  therefore the numerical evidence that we have gathered
is not yet conclusive regarding the robustness of the method.
C-CSP definitely deserves further analysis. 

The outline of this work is as follows. 
Section \ref{sec:preliminaries} is a quick overview on B-splines, NURBS, and isogeometric analysis.
In Section \ref{sec:method} we present a framework for isogeometric collocation and
the collocation schemes C-GP, LS-SP, C-ASP,  and the new C-CSP.
In Section \ref{sec:numerical-testing} we show some numerical tests of
C-CSP, focusing on the odd degree case, discuss its
robustness and compare it with the other collocation methods.
Finally, some conclusions and perspectives on future works are
detailed in Section \ref{sec:conclusion}.

\section{Preliminaries}\label{sec:preliminaries}

\subsection{B-splines}

Let us consider an interval $\hat{\Omega} \subset \mathbb{R}$. The B-splines
basis functions defined on $\hat{\Omega}$ are piecewise polynomials
that are built from a \emph{knot vector}, i.e.,  a vector with
non-decreasing entries $\Xi=[\xi_{1},\xi_{2}...,\xi_{n+p+1}]$,
where $n$ and $p$ are, respectively, the number of basis functions that
will be built from the knot vector and their polynomial degree.  
We name \emph{element} a \emph{knot span}  $(\xi_i, \xi_{i+1})$ having non-zero length, 
and we denote by $h$ the maximal length (the \emph{meshsize}).
A knot vector is said to be \textit{open} if its first and last
knot have multiplicity  $p+1$, i.e., each of them is repeated $p+1$ times.

Following \cite{DeBoor} and given a knot vector $\Xi$, univariate B-splines
basis functions $N_{i,p}$ are defined recursively as follows for $i=1,\ldots,n$: 
\begin{align}
& N_{i,0}(\xi)= 
\begin{cases}
1, &  \xi_{i}\leq \xi<\xi_{i+1},\\
0, &  \textrm{otherwise}, \\
\end{cases} 
\label{eq:Bsp}   \\ 
& N_{i,p}(\xi)=
\begin{cases}
  \dfrac{\xi-\xi_{i}}{\xi_{i+p}-\xi_{i}}N_{i,p-1}(\xi)+\dfrac{\xi_{i+p+1}-\xi}{\xi_{i+p+1}-\xi_{i+1}}N_{i+1,p-1}(\xi),
  & \xi_{i}\leq \xi<\xi_{i+p+1}, \\
0, & \textrm{otherwise,}
\end{cases} \nonumber
\end{align}
where we adopt the convention $\dfrac{0}{0}=0$;
note that the basis corresponding to an open knot vector
will be interpolatory in the first and last knot.

\begin{remark}
In this work we only consider  knot vectors whose internal
knots have multiplicity one: the associated B-splines/NURBS have then 
 global  $C^{p-1}$ regularity.
\end{remark}
We define by $\hat{S}^p=\textrm{span}\{N_{i,p} | i=1,...,n\}$
the space spanned by B-splines of degree $p$ and regularity $p-1$,  built
from a given knot vector $\Xi$. 
We also introduce the space of periodic B-splines, spanning the space
$\widetilde{S}^p=\{v\in \hat{S}^p | v(0)=v(1), \ v'(0)=v'(1),..., \ v^{(p-1)}(0)=v^{(p-1)}(1)\}$;
interestingly, the dimension of $\widetilde{S}^p$ equals
the number of elements of the underlying knot vector $\Xi$, a property that will come in handy
later on.

Multivariate splines spaces can be constructed from univariate spaces by means of tensor products. 
For example, a B-splines space in two dimensions can be defined by considering the knot vectors
$\Xi=[\xi_{1},\xi_{2}...,\xi_{n+p+1}]$ and $\Lambda=[\eta_{1},\eta_{2},...,\eta_{m+q+1}]$,
and defining
$\hat{S}^{p,q}=\textrm{span}\{N_{i,p}(\xi)M_{j,q}(\eta),  i=1,...,n, j=1,...,m  \}$.
In the following, it will be useful
to refer to the basis functions spanning $\hat{S}^{p,q}$ with a single running index $k$
ranging from $1$ to $n \times m$, i.e.
\begin{equation}\label{eq:S_bidim_one_idx}
\hat{S}^{p,q} = \textrm{span}\{\varphi_{k}^{p,q}(\xi,\eta)=N_{i,p}(\xi)M_{j,q}(\eta)\, |\, k = i+(j-1)m, \, i=1,\ldots,n, \, j=1,\ldots,m  \}. 
\end{equation}

\subsection{NURBS}

Non-uniform rational B-splines (NURBS, cf. \cite{Piegl:2012}) 
are defined for the purpose of describing geometries of practical interest like conic sections,
{see e.g. Problem \ref{prob:dirichlet} in next section}.  
The definition of a generic bivariate NURBS function on the parametric
square $\hat \Omega$ is
\[
\forall(\xi, \eta) \in \hat{\Omega},\quad
 R_{i,j}^{p,q}(\xi,\eta)=
 \dfrac{N_{i,p}(\xi)M_{j,q}(\eta)w_{i,j}}{\sum_{\widehat{i}=1}^{n}\sum_{\widehat{j}=1}^{m}N_{\widehat{i},p}(\xi)M_{\widehat{j},q}(\eta)w_{\hat{i},\hat{j}}} 
 \forall i=1,...,n, \forall j=1,...,m
\]
where $w_{i}$ are suitable weights, and $N_{i,p}(\xi), M_{j,q}(\eta)$ are the univariate 
B-splines basis functions defined in \eqref{eq:Bsp}.
Similarly to \eqref{eq:S_bidim_one_idx} we also introduce a single running index $k=1,\ldots,n \times m$
to refer to the NURBS basis, i.e.,
\[
 R_{k}^{p,q}(\xi,\eta)= R_{i,j}^{p,q}(\xi,\eta), \text{ with }  k = i+(j-1)m, \, i=1,\ldots,n, \, j=1,\ldots,m. 
\]

\section{Isogeometric collocation and the choice of the collocation points}
\label{sec:method}

\subsection{Isogeometric collocation}

In our numerical tests we will consider both one-dimensional and two-dimensional
elliptic problems, which we now introduce.

\begin{prob}[One-dimensional Dirichlet boundary problem]\label{prob:dirichlet1d}
  Find $u:[0,1]\rightarrow \mathbb{R}$ such that
  \begin{equation}
    \begin{cases}
      u''(x)+a_{1}(x)u'(x)+a_{0}(x)u(x)=f(x) \ \ \ \forall x \in (0,1)\\
      u(0)=u(1)=0 
    \end{cases}
    \label{dir1d}
  \end{equation}
  where $a_0, a_1, f:[0,1]\rightarrow \mathbb{R}$ are sufficiently regular functions.
\end{prob}
We assume that this problem has a unique smooth solution.
We then look for an approximate solution
$u_{h}(x)=\sum_{i=1}^{n}c_iN_{i,p}(x) \in \hat{S}^p$,
that complies with the boundary conditions $u(0)=u(1)=0$ (i.e. $c_1=c_n=0$, 
given the interpolatory property of open knot vectors at the first and last knot), 
and that satisfies  \eqref{dir1d} in $n-2$ \textit{collocation points}
$\{\tau_1, ..., \tau_{n-2}\}$ that need to be specified, i.e.
\begin{equation}\label{eq:1D-collocation}
u_{h}''(\tau_i)+a_{1}u_{h}'(\tau_i)+a_{0}u_{h}(\tau_i)=f(\tau_i),  \qquad \forall i=1,...,n-2.
\end{equation}
The coefficients $c_2, \ldots, c_{n-1}$ are then computed by solving the linear system
obtained by inserting the expansion $u_{h}(x)=\sum_{i=1}^{n}c_iN_{i,p}(x)$ into \eqref{eq:1D-collocation}.
We also shall introduce a periodic version of Problem \ref{prob:dirichlet1d}, 
which we consider because it is particularly simple
to set up a collocation scheme for it, due to the already-mentioned fact
that the number of degrees-of-freedom $n$ of $\widetilde{S}^p$ 
(hence the number of collocation points to be used) is identical to
the number of elements of $\Xi$.
\begin{prob}[One-dimensional periodic boundary problem]\label{prob:periodic}
  Find $u: \mathbb{R}\rightarrow \mathbb{R}$ such that
  \begin{equation}
    \begin{cases}
      u''(x)+a_{1}u'(x)+a_{0}u(x)=f(x) 	& \forall x \in \mathbb{R},\\
      u(x)=u(1+x) 						&  \forall x \in \mathbb{R} ,
    \end{cases}
    \label{per1d}
  \end{equation}
  where $a_0, a_1$ and $f$ are sufficiently regular periodic functions.
\end{prob}
We assume again that this problem has a unique (periodic) smooth solution.
{Note that the periodic problem is not well-posed if $a_0$ is null.}
The B-splines approximation of the solution of 
\eqref{per1d} is therefore $u_{h}\in \widetilde{S}^p$ such that 
\begin{equation}\label{eq:1D-collocation-per}
u_{h}''(\tau_i)+a_{1}u_{h}'(\tau_i)+a_{0}u_{h}(\tau_i)=f(\tau_i)  \qquad \forall i=1,...,n,
\end{equation}
for suitably chosen collocation points $\{\tau_1, ..., \tau_{n}\}$
with periodic distribution on $[0,1]$.

Finally, we also consider the two-dimensional Poisson equation, 
that we will solve by a multivariate collocation scheme constructed 
by tensorizing univariate sets of collocation points. More specifically, we denote by 
$\Omega\subset \mathbb{R}^2$
a domain described by a NURBS parametrization $\mathbf{F}:\hat{\Omega}\rightarrow \Omega$,
where $\hat{\Omega}=[0,1]\times[0,1]$ and 
\[
\mathbf{F}(\xi,\eta) = \sum_{k=1}^{n \times m} \mathbf{P}_k R_{k}^{p,q}(\xi,\eta), \quad \mathbf{P}_k \in \mathbb{R}^2,
\]
we let $\Gamma$ denote the boundary of $\Omega$, and we consider the Dirichlet problem
\begin{prob}[Two-dimensional Dirichlet boundary problem]\label{prob:dirichlet}
Find $u:\Omega \rightarrow \mathbb{R}$ such that 
\begin{equation}
\begin{cases}
-\Delta u=f & \textrm{in} \ \Omega,\\
u=0   & \textrm{on} \ \Gamma,
\end{cases}
\label{dir2d}
\end{equation}
where $f:\Omega\rightarrow \mathbb{R}$ is a sufficiently regular function.   
\end{prob}
Again, we assume that this problem has a unique smooth solution.
Following the isogeometric paradigm, the discrete solution $u_{h}$ is sought in the isogeometric space
\[
u_{h} \in  S^{p,q} = 
\text{span} \left \{ R_{k}^{p,q} \circ \mathbf{F}^{-1}, \, \forall k = i+(j-1)m, \, i=1,...,n, \, j=1,...,m \right \}
\]
cf. \eqref{eq:S_bidim_one_idx}, and the collocation points are the image 
through $\mathbf{F}(\cdot)$ of a tensor-product grid of
collocation points on $[0,1]^2$. 
The collocation method is then obtained as for the univariate case.

\subsection{Greville points and C-GP}

\begin{figure}[tp]
 \centering
 \includegraphics[width=.8\textwidth]{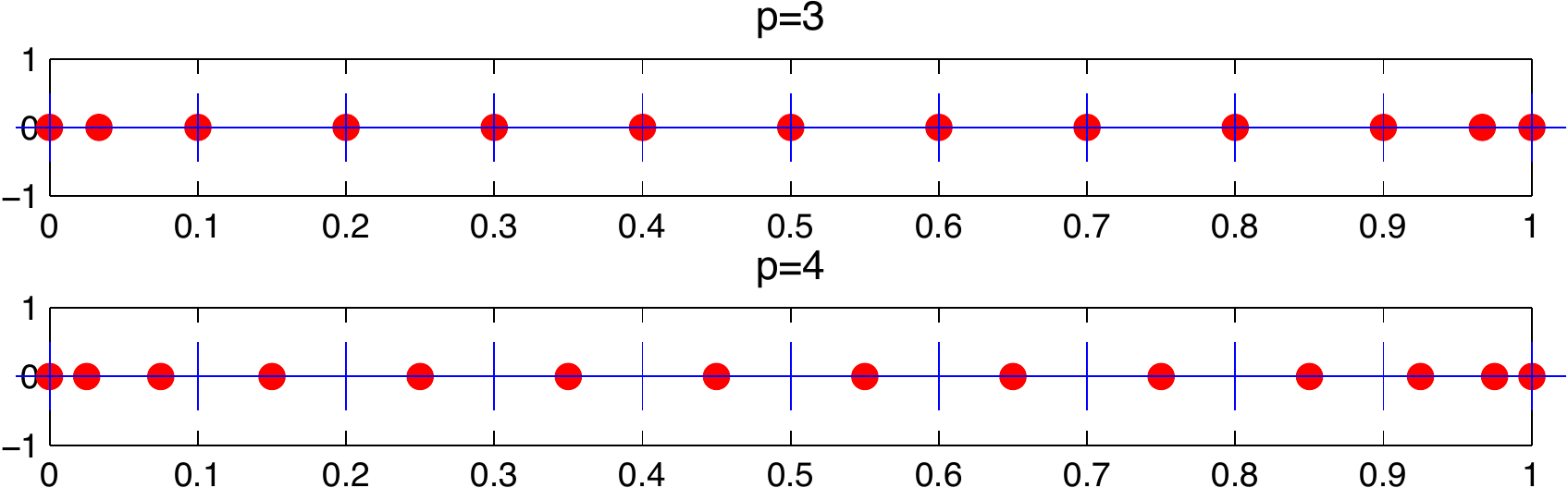}
 \caption{Examples of Greville points computed from an open  knot vector:
   p=3 and p=4. The interior Greville points are used as collocation
   points in the C-GP scheme.}
 \label{ex_grev}
\end{figure}

Greville points, or abscissas, for $p$-degree  B-splines associated to
a knot vector $\Xi=\{\xi_1,...,\xi_{n+p+1}\}$ are defined as 
\[%\begin{equation}
{\tau^{GP}_i}=\frac{\xi_{i+1}+...+\xi_{i+p}}{p}, \quad \forall i=1,...,n,
\]%\end{equation}
see Figure $\ref{ex_grev}$ for an example computed from a open uniform
knot vector and degree $p=3$ and $p=4$. 
For an open knot vector the first and last Greville point coincide
with the first and last knot $\xi_1$ and $\xi_{n+p+1} $. A common
collocation scheme for second-order
boundary value problems,  as proposed in \cite{Auricchio:2010},  uses
as collocation points the $n-2$ internal Greville points. For brevity,  
this is denoted Collocation at Greville Points, C-GP.

\begin{table}
  \begin{center}
    \begin{tabular}{c|c|cc|cc|c}
      & Galerkin & \multicolumn{2}{c}{C-GP} & \multicolumn{2}{c}{LS-SP
        and C-CSP} & C-ASP\\
      & & Odd $p$ & Even $p$ & Odd $p$& Even $p$&  \\
      \hline
      $L^2 $ & $p+1$ &$p-1$ & $p$&$p+1 $ & $p$& $p$\\
      $H^1 $ & $p$ &$p-1$ & $p$&$p$ & $p$& $p$\\
      $H^2$ & $p-1$ &$p-1$ & $p-1$& $p-1$ & $p-1$ & $p-1$\\
      \hline
    \end{tabular}
  \end{center}
  \caption{Comparisons of orders of convergence: Galerkin, C-GP,
    LS-SP,  C-CSP and C-ASP.}
  \label{compar}
\end{table}

In Table~\ref{compar} we report the orders of
convergence of C-GP and of the other methods considered in this paper. 
The convergence rate of C-GP in $L^2$ norm
is $p-1$ when odd degree are used and $p$ when even degree
B-Splines are used as discussed earlier (i.e. two-orders and one-order suboptimal, respectively). 
The error in $H^1$ norm converges with the same orders of 
the $L^2$ norm, and is therefore optimal for even degrees and one-order suboptimal for odd degrees. 
The error measured in $H^2$ norm is instead optimal for every degree.

\subsection{Cauchy-Galerkin points and superconvergent points for the second derivative of the Galerkin solution}\label{sec:sppoints}

Following \cite{gomez2016variational} and \cite{anitescu2015isogeometric},
we now introduce the Cauchy-Galerkin points  and the second-derivative superconvergent points 
for the Galerkin solution of Problem \ref{prob:dirichlet1d}, which will be used to construct
a  collocation or least-squares method. Assume for a moment
that $a_0=a_1=0$ in Problem \ref{prob:dirichlet1d}, i.e., consider 
\begin{equation}
  \label{eq:problema-per-grafico-superconvergenza}
  \begin{cases}
    - u''(x)=f(x) & \ \ \ \forall x \in (0,1)\\
    u(0)=u(1)=0,
\end{cases}
 \end{equation}
and let $u_h^*$ be the approximated solution given by the Galerkin method based on B-splines. 
The Cauchy-Galerkin points are  collocation points where the
Galerkin residual,  in this case  $D^2(u-u_h^*)$, is
zero. Since these points are unknown a-priori, one can look for a
high-order approximation of them, i.e, the so-called superconvergent points. In general, the points
$\Psi_h=\{\psi_{h,1}, ..., \psi_{h,w} \}$ with $w \in \mathbb{N}, w>0,$ are said to be
superconvergent points for the $j$-th derivative of $u$ if
\begin{equation}
\Big[\sum_{\psi_{h,i} \in \Psi_h}\big[D^j(u-u_h^*)(\psi_{h,i})\big]^2\Big]^{\frac{1}{2}}\leq Ch^{p+1-j+k},
\quad \forall i=1\ldots,w,
\label{sp}
\end{equation}
where $k>0$, $j\geq 0$, $C$ is a constant, $h$ is the meshsize of the knot vector, and 
$p$ is the degree of the B-splines. Here we are therefore
interested in the case $k=1$ and $j=2$. 

\begin{table}
\begin{center}
\begin{tabular}{c|c}
Degree & Second derivative SP \\
\hline
p=3 & $\frac{-1}{\sqrt{3}}, \frac{1}{\sqrt{3}}$\\
p=4 & -1,0,1\\
p=5 & $\pm\frac{\sqrt{225-30\sqrt{30}}}{15}$\\
p=6 & -1,0,1\\
p=7 & $\pm 0.504918567512$\\
\hline
\end{tabular}
\end{center}
\caption{On the reference element $[-1,1]$, location of superconvergent
  points for the second derivative (from \cite{gomez2016variational}).}
\label{sptable}
\end{table}

Finding the location of the superconvergent points
is in general an open problem as well.
Under the assumption that the superconvergent points are element-invariant 
(i.e., images by affine mapping of points on a reference element)  
their locations have been estimated in \cite{gomez2016variational} and
are reported in Table \ref{sptable} for a reference element
$[-1,1]$. The same  points are estimated in
\cite{anitescu2015isogeometric} under a similar periodicity
assumption. Both assumptions do not hold true
in many cases of interest.  {An alternative superconvergence theory can be found in 
\cite{wahlbin1995superconvergence}, based on a mesh symmetry
assumption; this hypothesis however} does not hold true for elements close to the boundary.

\begin{figure}[tp]
\centering
\includegraphics[width=0.7\textwidth]{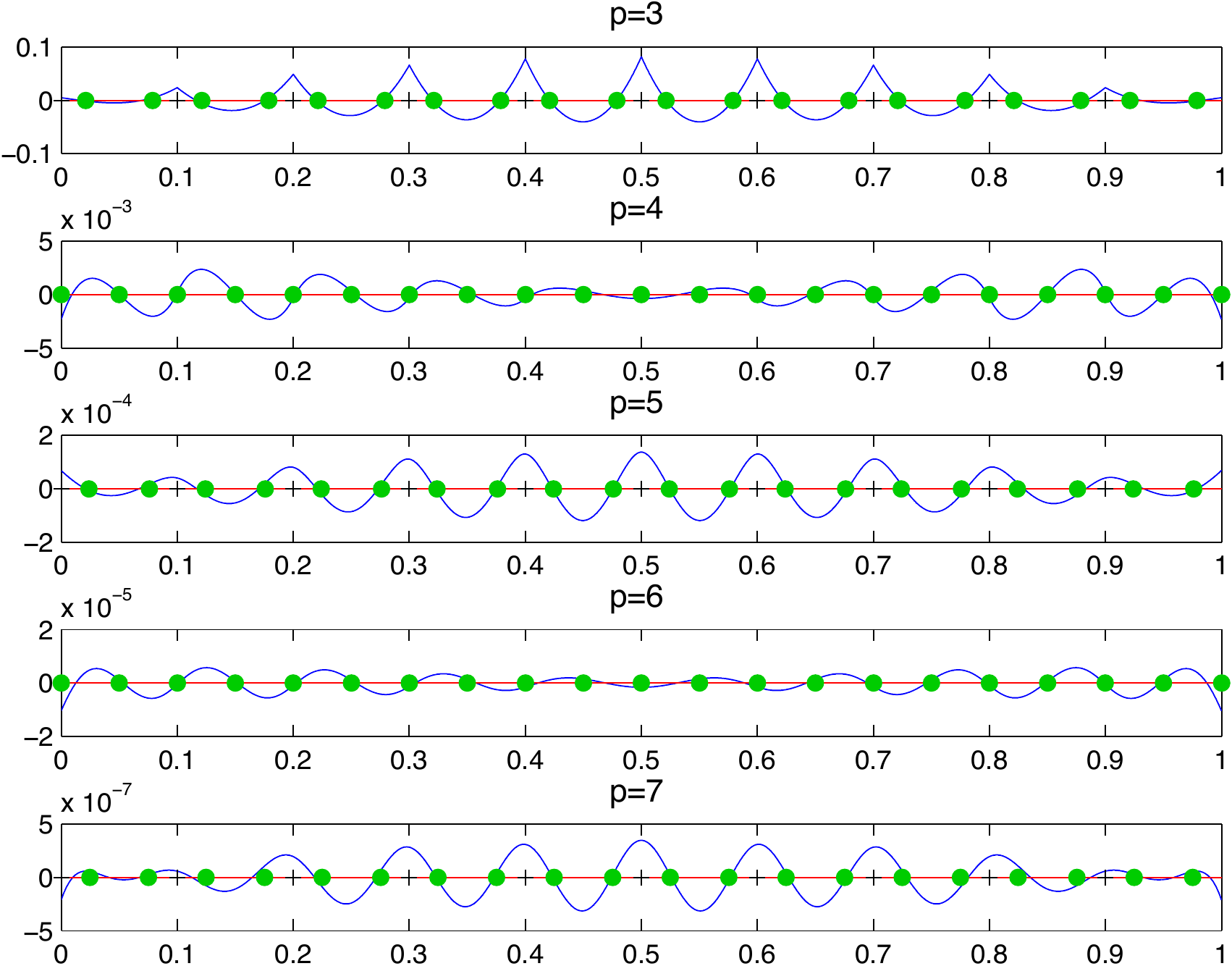}
 \caption{Plot of   $D^2(u-u_h^*)$, equivalent to the residual of problem
   \eqref{eq:problema-per-grafico-superconvergenza},  and ``surrogate''
   superconvergent points (green dots), on a mesh with 10 elements.}
 \label{res_dir_10}
\end{figure}

\begin{figure}[tp]
\centering
\includegraphics[width=0.7\textwidth]{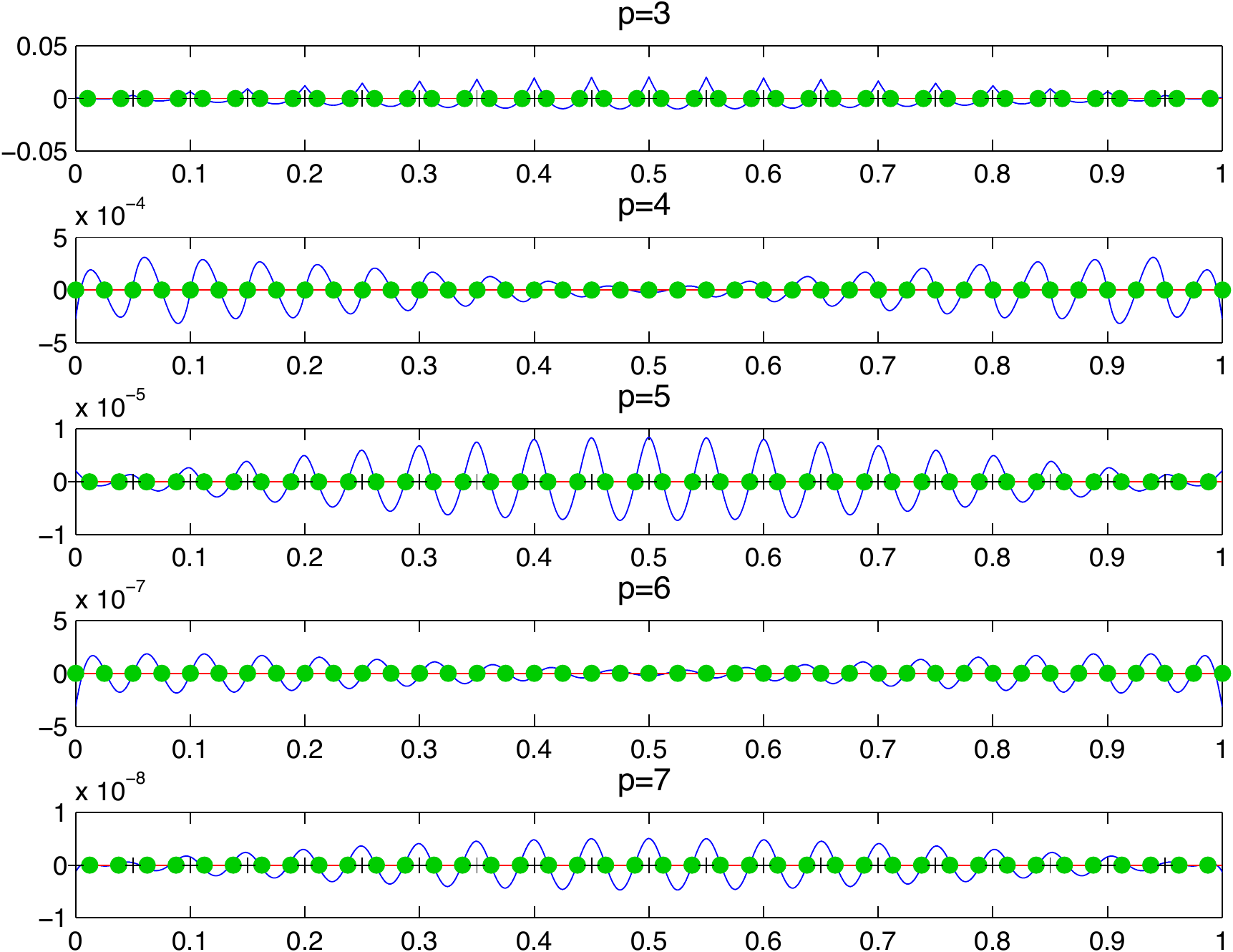}
 \caption{Plot of   $D^2(u-u_h^*)$, equivalent to the residual of problem
   \eqref{eq:problema-per-grafico-superconvergenza},  and ``surrogate''
   superconvergent points (green dots), on a mesh with 20 elements.}
 \label{res_dir_20}
\end{figure}

\begin{figure}[tp]
\centering
\includegraphics[width=0.7\textwidth]{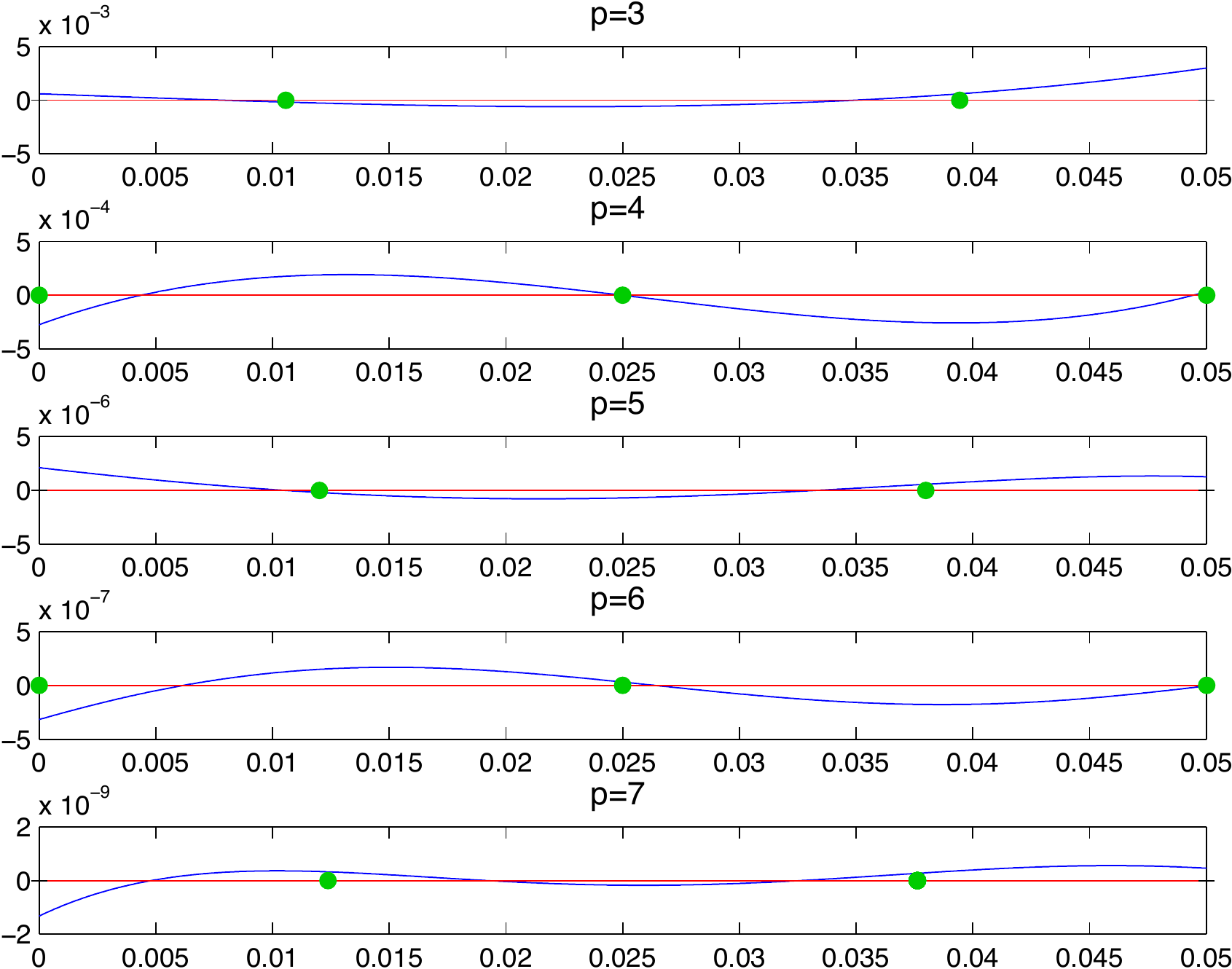}
\caption{Zoom on the first element of Figure \ref{res_dir_20}.}
\label{res_dir_20_zoom}
\end{figure}

\begin{figure}[tp]
\centering
\includegraphics[width=0.7\textwidth]{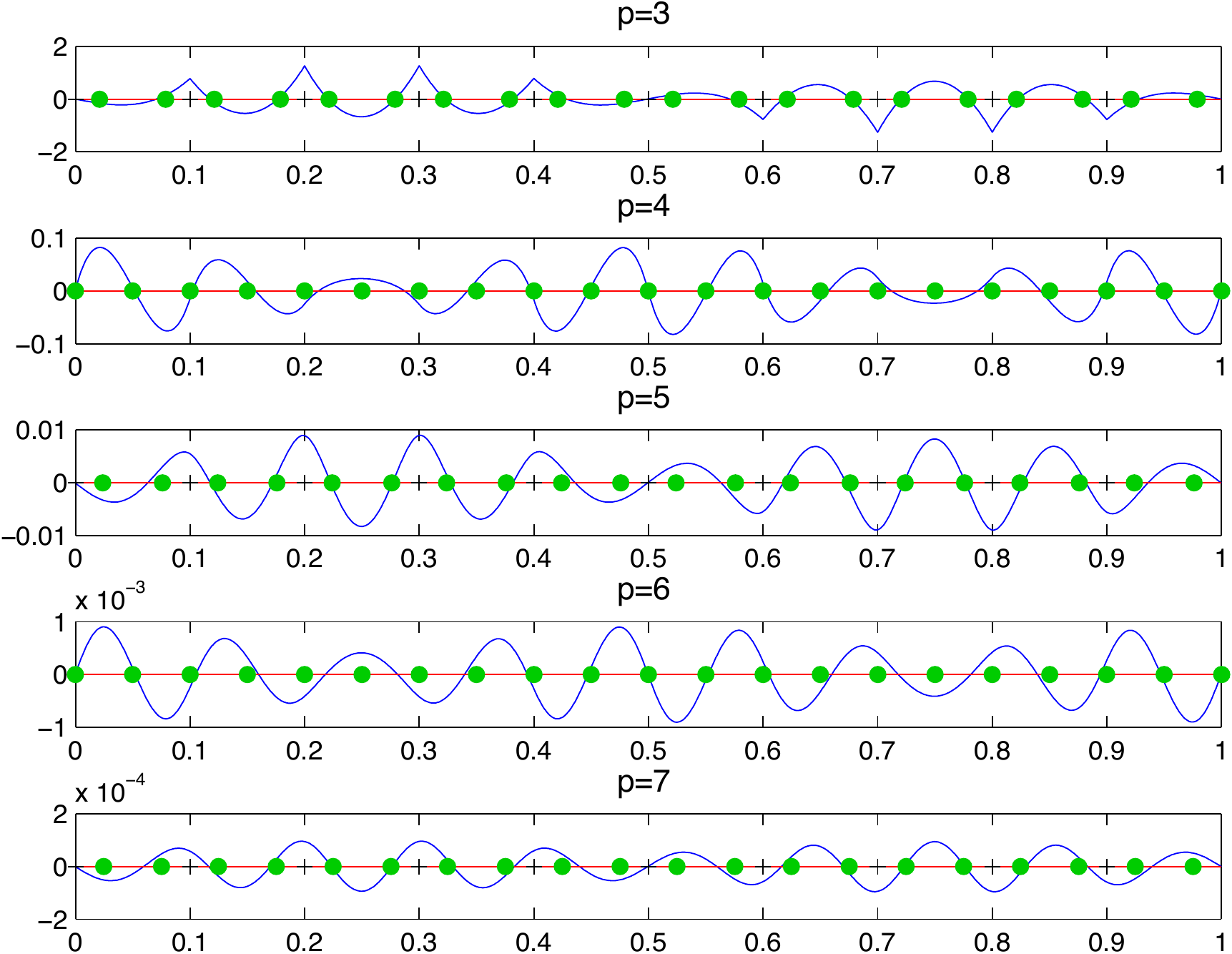}
 \caption{Residuals of periodic problem with 10 elements.}
 \label{res_per_10}
\end{figure}

Following \cite{anitescu2015isogeometric}  and
\cite{gomez2016variational}, since  we do not have access to the ``true'' superconvergent
points, we use  the points in  Table  $\ref{sptable}$, linearly mapped to the
generic element, as  ``surrogate'' superconvergent points in one-dimension. 
How well do these ``surrogate'' superconvergent points
approximate the Cauchy-Galerkin points? This is the main question and
some qualitative answer {can be found} in Figures \ref{res_dir_10} and
\ref{res_dir_20}, that show $D^2(u-u_h^*)$ 
for equation \eqref{eq:problema-per-grafico-superconvergenza} with
$f(x)=\sin(\pi x)$, over a mesh with $10$ and $20$ elements and $p=3,\ldots,7$, as well as
the ``surrogate'' superconvergent points for each degree of approximation:
for odd degrees, a non-negligible discrepancy is evident at the 
{boundaries of the domain}, 
and for even degrees this occurs also at the middle of the
interval. Figure \ref{res_dir_20_zoom} is a zoom of the first element in Figure \ref{res_dir_20}. 
For completeness, Figure \ref{res_per_10} 
shows the residual for the Periodic Problem \ref{prob:periodic} with
$a_0=a_1=1$, and $f(x)=(1+4\pi^2)\sin(\pi x)+2\pi\cos(2\pi x)$   
over a mesh with 10  elements and $p=3,\ldots,7$, as well
as  the ``surrogate'' superconvergent points for each degree of approximation.  In this
case, the mismatch between the zeros of the residual and the
``surrogate''  superconvergent points is higher in correspondence of a
smaller residual. Note that the residual is not periodic at the
element scale.

For easiness of exposition, from now on we refer to  the ``surrogate''
superconvergent points  simply as superconvergent points, although this might not be technically
true.

For multi-dimensional problems on a NURBS single-patch geometry, the
superconvergent points
can be obtained by further mapping the tensor
product of one-dimensional superconvergent points through the geometry map $\mathbf{F}$ in
the physical domain.  Clearly, the same considerations of 
the one-dimensional case are valid.

\subsection{Least-Squares at Superconvergent Points (LS-SP)}

As already mentioned, the Least-Squares at Superconvergent Points method (LS-SP) has
been introduced by \cite{anitescu2015isogeometric}.  In this method
all the superconvergent points are used as
collocation points. As it can be seen in Table \ref{sptable}, there are at least two
superconvergent points per element; if we take all of them as
collocation points, we obtain an overdetermined system of
equations if the number of elements is large enough: 
such linear system is then solved in a least-squares sense, leading
to a method which is not strictly a collocation method.
The order of convergence of the method as measured in numerical tests is reported in Table
\ref{compar}: note that it is optimal for odd degrees and one-order sub-optimal
in $L^2$ for even degrees, while it is optimal regardless of the parity of $p$ in $H^1$ and $H^2$ norm. 

\begin{figure}[tp]
\begin{center}
\includegraphics[height=0.35\textheight]{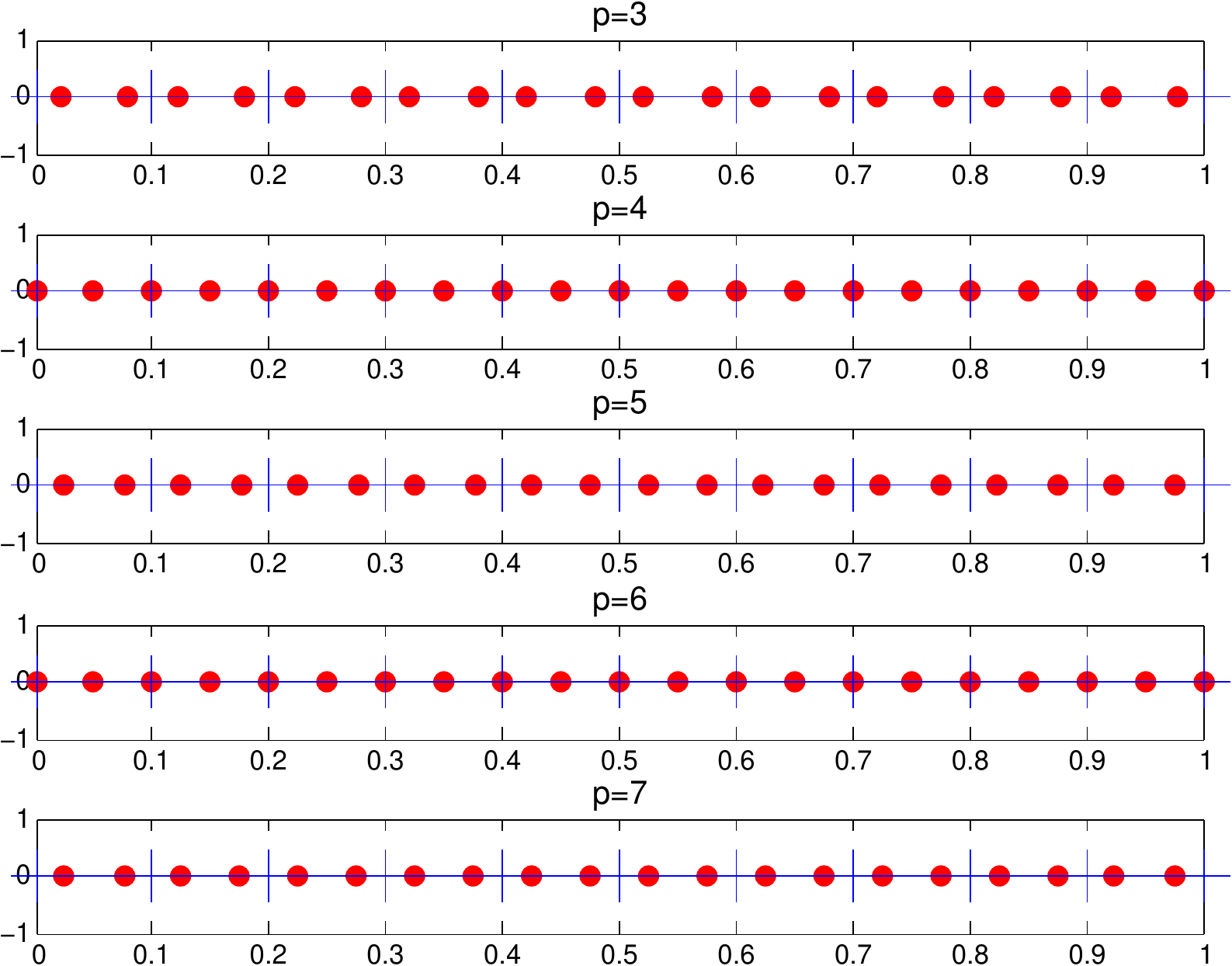}
\caption{Superconvergent points for  $p=3,\ldots,7$ on a knot vector with 10 elements.}
\label{SP}
\end{center}
\end{figure}

Figure $\ref{SP}$ shows the superconvergent points for $p=3,\ldots,7$ 
on a knot vector with 10 elements. Observe that the same least-squares formulation
can accommodate for both Dirichlet problems (i.e., open knot vectors) and periodic problems.

\subsection{Collocation at Alternating Superconvergent Points (C-ASP)}

\begin{figure}[tp]
\begin{center}
\subfloat[][Case $p$ odd. \label{aspperiodicodd}]
{\includegraphics[height=0.02\textheight]{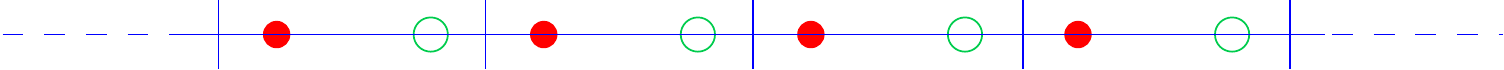}}\quad
\subfloat[][Case $p$ even.\label{aspperiodiceven}]{\includegraphics[height=0.02\textheight]{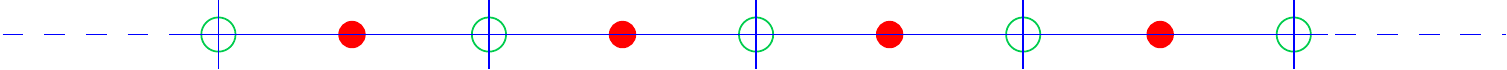}}
\caption{Example of C-ASP points for the periodic Problem \ref{prob:periodic}. 
  The collocation points are marked with full red dots, while the remaining superconvergent points are 
  displayed with green circles. In this case, C-ASP and C-GP coincide for even degrees.}
\label{ACPpoints-periodic}
\end{center}
\end{figure}

\begin{figure}[tp]
\begin{center}
\includegraphics[height=0.35\textheight]{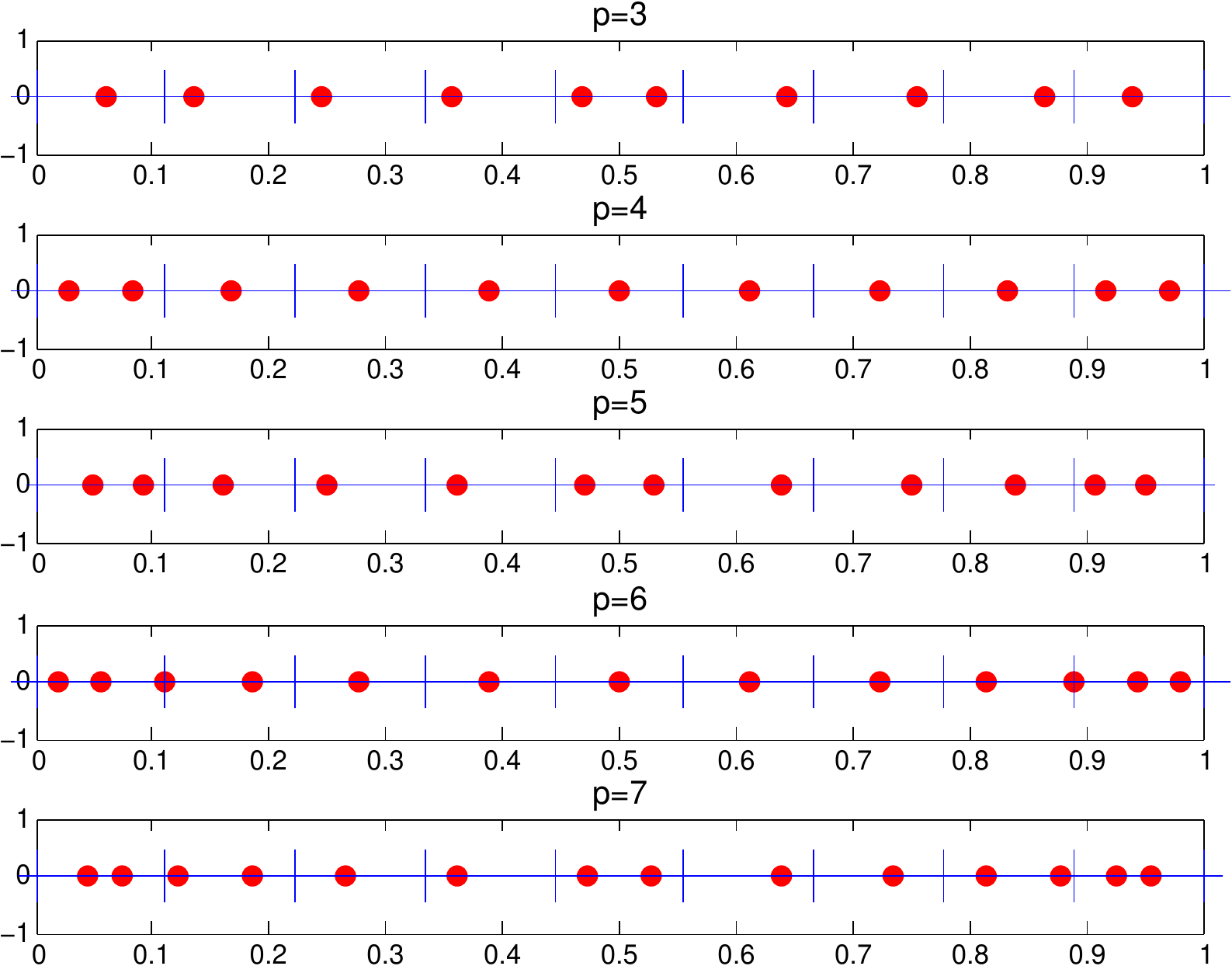}
\caption{Example of C-ASP points for the Dirichlet problem \ref{prob:dirichlet1d} over a knot vector with 9 elements.
The points adjacent to the boundary are obtained according to Algorithm 1 of \cite{gomez2016variational}.}
\label{ACPpoints}
\end{center}
\end{figure}

C-ASP is a collocation method introduced in \cite{gomez2016variational}. 
In this method, a subset of superconvergent points with cardinality
equal to the number of degrees-of-freedom is employed
as set of collocation points. The authors of \cite{gomez2016variational} select 
a subset of the superconvergent points in such a way that 
every element of the knot span contains at least one collocation point;
note that this roughly means considering every other superconvergent
point, hence the name we give to the method.
Because we need to select as many collocation points as degrees of freedom, the easiest case is when one considers
the periodic Problem \ref{prob:periodic}, for which the number of elements is identical to the number of degrees-of-freedom,
so that exactly one superconvergent point per element is selected, see
Figure \ref{ACPpoints-periodic} (this case is not considered in
\cite{gomez2016variational}). Note that for even $p$ one possibility is
then to select the midpoint of each element, i.e., the Greville
points for the uniform knot vector, see Figure \ref{aspperiodiceven}. 
For the Dirichlet Problem \ref{prob:dirichlet1d},  
one needs instead to select $n_{el} +p -2$ collocation points on a mesh of $n_{el} $ elements. 
To this end, an ad-hoc algorithm is presented in \cite{gomez2016variational}
that selects suitable superconvergence points in the internal part of the domain, 
and ``blends them'' with Greville points on the elements close to the boundary, as
can be seen in Figure \ref{ACPpoints}. Note that other choices for the elements
close to the boundary can be envisaged, which however do not affect the convergence
order of the method, see \cite{montardini:MSthesis}. 

The convergence orders of C-ASP assessed numerically by \cite{gomez2016variational}
are also reported in Table \ref{compar}. 
Note in particular that the $L^{2}$ order of convergence for C-ASP is $p$ 
regardless of the parity of $p$, i.e., one-order suboptimal, 
while the $H^1$ and $H^2$ orders of convergence are optimal, again regardless of the parity of $p$.

\subsection{Collocation on Clustered Superconvergent Points (C-CSP)}

We now describe a new choice of collocation points among the
superconvergent points, alternative to C-ASP, 
which we name Collocation on Clustered Superconvergent Points
(C-CSP). 

To understand our approach, we describe it first in the simplest setting, i.e., 
the periodic Problem \ref{prob:periodic} with even number of elements and odd degree $p$. 
We look for a periodic distribution of collocation points 
{which is furthermore symmetric} at the element scale.
This can be achieved by selecting two superconvergent points in an element and then skipping the
following one, as depicted in Figure \ref{periodic_ccsp_odd}.
Surprisingly, the order of convergence of C-CSP in this
case is optimal, cf. the numerical results in Section \ref{sec:numerical-testing}, 
Figure \ref{CSPperiodicerror}.
For even degrees, we have experimented different selections  of sets of 
superconvergent points, preserving  periodicity and some local
symmetry, two of which are depicted in Figure \ref{periodic_ccsp_even} 
(observe that with the first one we end up with Greville points again). 
In all cases, we have measured numerically one-order suboptimal convergence in $L^2$, i.e., 
we do not see improvements with respect to C-GP, LS-SP or  C-ASP, 
see the numerical results in Section \ref{sec:numerical-testing}, 
Figure \ref{CSPperiodicerror-even}.
At this point, only the odd-degree C-CSP seems to deserve further
interest, and we will restrict to this case in the remaining of the
paper. How to use efficiently the superconvergent points in an even
degree splines collocation scheme remains an open problem.

The next step is to extend (odd-degree) C-CSP to the open knot vector
to solve the Dirichlet Problem \ref{prob:dirichlet1d}.  To this
end, we need to include additional points, which are taken among the
other superconvergent points populating the elements close to the
boundary, and trying to preserve symmetry, see Figures \ref{CSPodd}
and \ref{CSPeven}.  Note that {when the number of elements is even 
the procedure just described will not yield a globally symmetric
distribution of collocation points, cf. Figure \ref{CSPeven}.
We can however restore symmetry with a little modification of the collocation
approach: we add one (or a few) points to the collocation set to
restore symmetry of the collocation scheme, and average the equations
corresponding to the points located at the center of the domain in
order to match the number of unknown.  This procedure is depicted in
Figure \ref{CSP_sym} for $p=3$.

The order of convergence of C-CSP on regular meshes, reported in
Table \ref{compar}, is the same of LS-SP, i.e., optimal for odd
degrees. As already mentioned, all our attempts to extend it to even degree splines has
produced one order suboptimal convergence.
These convergence rates have been measured by running the 
numerical benchmarks detailed in Section \ref{sec:numerical-testing}, 
also covering the symmetric variant of Figure \ref{CSPeven}.

\begin{figure}[tp]
\begin{center}
{\includegraphics[height=0.02\textheight]{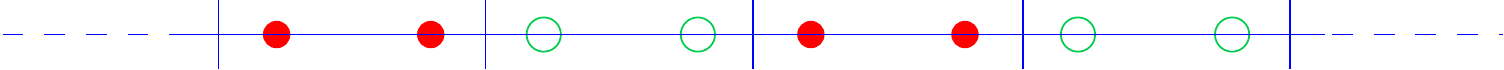}}
\caption{Periodic C-CSP stencil for odd degree: the collocation points are
  marked with full red dots, while the remaining superconvergent
  points are  displayed with green circles.}
\label{periodic_ccsp_odd}
\end{center}
\end{figure}

\begin{figure}[tp]
\begin{center}
{\includegraphics[height=0.02\textheight]{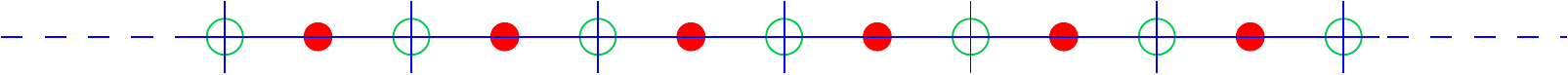}}\\
{\includegraphics[height=0.02\textheight]{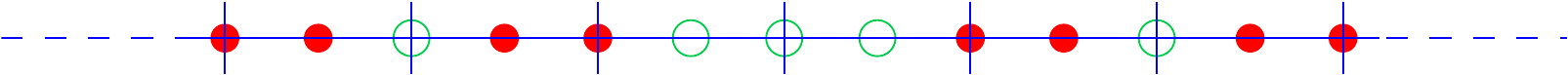}}
\caption{Attempts of C-CSP stencil for even degree: the collocation points are marked with full red dots, 
while the remaining superconvergent points are 
displayed with green circles. The construction at the top leads to Greville points, while the one at the bottom yields symmetry
at a macro-element level.}
\label{periodic_ccsp_even}
\end{center}
\end{figure}

\begin{figure}[tp]
\begin{center}
\includegraphics[height=0.35\textheight]{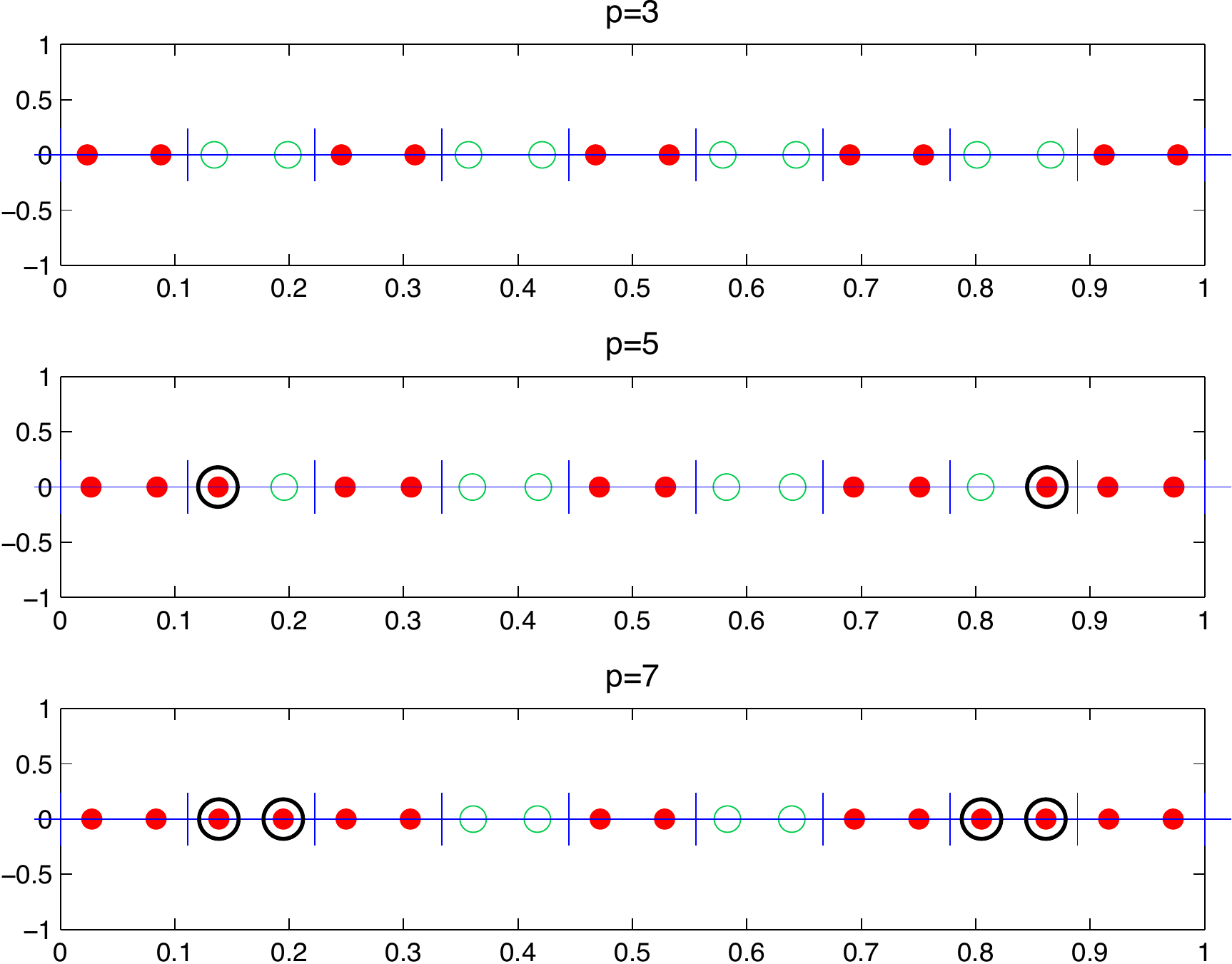}
\caption{C-CSP points for a Dirichlet problem solved on a mesh with 9 elements 
(\emph{odd number of elements}, {leading to a symmetric set of point}):
the collocation points are marked with full red dots, while the remaining superconvergent points are 
displayed with green circles. Black dots represent the points added with respect to the periodic stencil.
 }
\label{CSPodd}
\end{center}
\end{figure}

\begin{figure}[tp]
\begin{center}
\includegraphics[height=0.35\textheight]{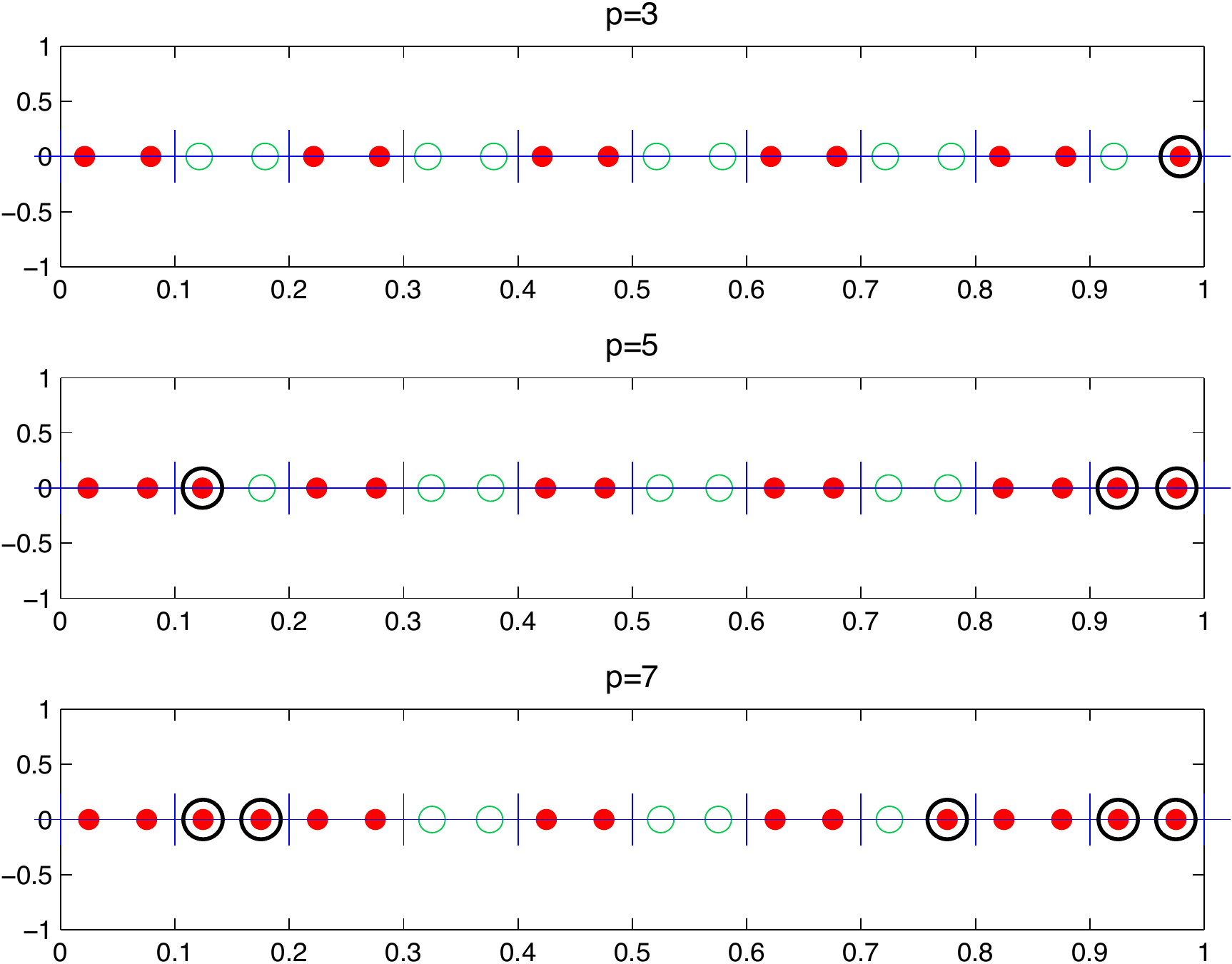}
\caption{C-CSP points for a Dirichlet problem solved on a mesh with 10 elements 
 (\emph{even number of elements}, {leading to a non-symmetric set of point}):
the collocation points are marked with full red dots, while the remaining superconvergent points are 
displayed with green circles. Black dots represent the points added with respect to the periodic stencil.
}
\label{CSPeven}
\end{center}
\end{figure}

\begin{figure}[tp]
\begin{center}
\includegraphics[height=0.35\textheight]{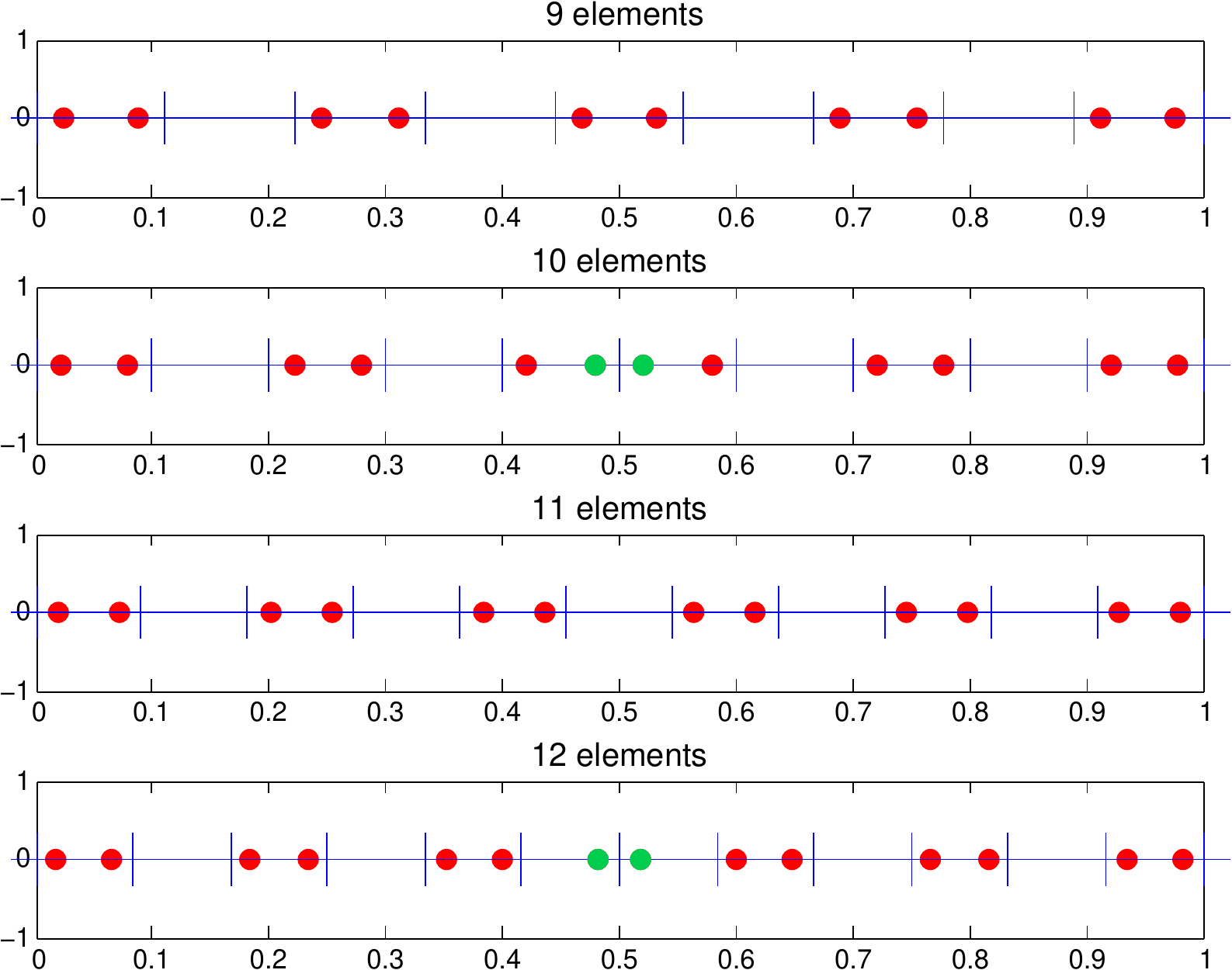}
\caption{C-CSP symmetric-variant points for a Dirichlet problem with $p=3$:
the collocation points are marked with  red dots, while the superconvergent points whose equations have to be averaged are 
displayed with green dots.}
\label{CSP_sym}
\end{center}
\end{figure}

\section{Numerical tests}
\label{sec:numerical-testing}

This section is devoted to the numerical benchmarking of the new C-CSP
method, and its comparison to the other approaches recalled in Section \ref{sec:method}. 
For conciseness, we do not show convergence results in $L^\infty$ norm, which we found to 
be identical to the ones in $L^2$ norm in each of the tests reported
below, {see \cite{montardini:MSthesis}}.

We begin by testing C-CSP on the periodic Problem \ref{prob:periodic}, with
$a_0=a_1=1$ and with $f(x)=(1+4\pi^2)\sin(2\pi x)+2\pi \cos(2\pi x)$, 
whose solution is $u(x)=\sin(2\pi x)$. 
As previously discussed, this is the only test for which we present results for 
even degrees $p$: we see from the plots in Figures  \ref{CSPperiodicerror} and \ref{CSPperiodicerror-even} that
the orders of convergence for the $L^{2}$ norm of the error are optimal,  
i.e. equal to $p+1$, for odd values of $p$, while for even $p$ the 
measured convergence rate is only $p$, i.e. one-order suboptimal.

\begin{figure}[tp]
\begin{center}
\includegraphics[width=0.44\linewidth]{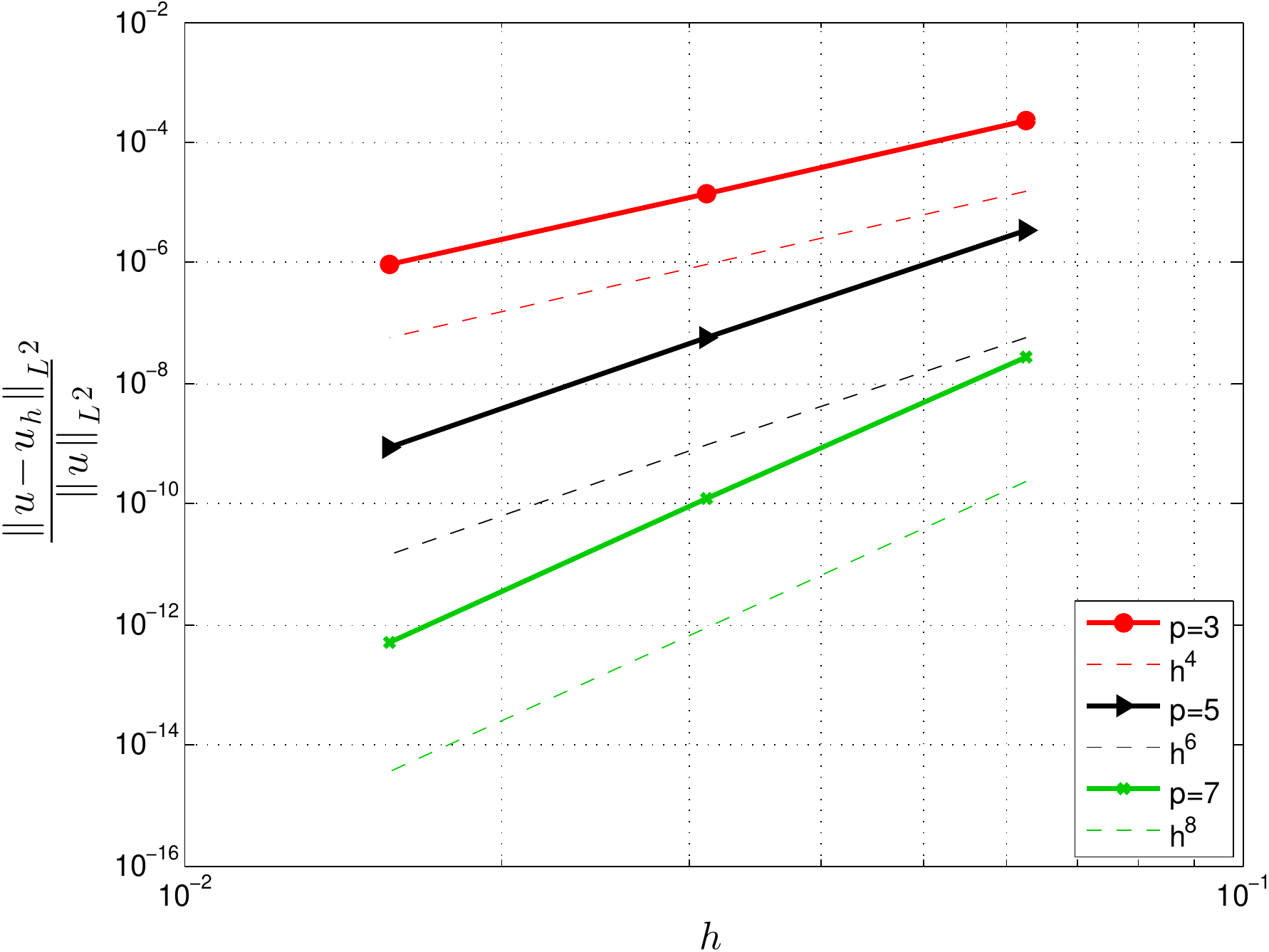}
\includegraphics[width=0.44\linewidth]{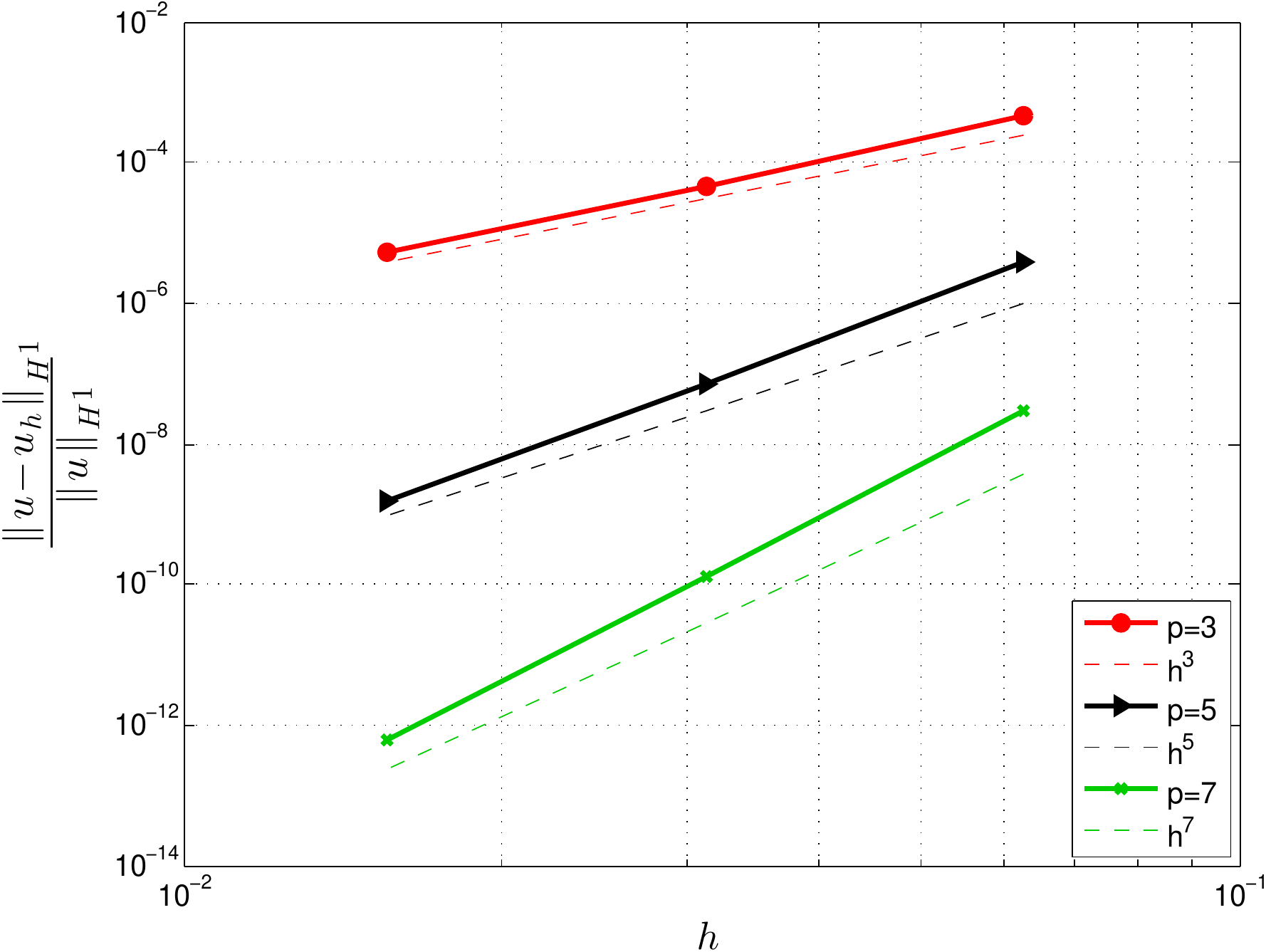}
\caption{$L^{2}$ and $H^{1}$  error plot: C-CSP periodic problem (odd $p$).}
\label{CSPperiodicerror}
\end{center}
\end{figure}

\begin{figure}[tp]
\begin{center}
\includegraphics[width=0.44\linewidth]{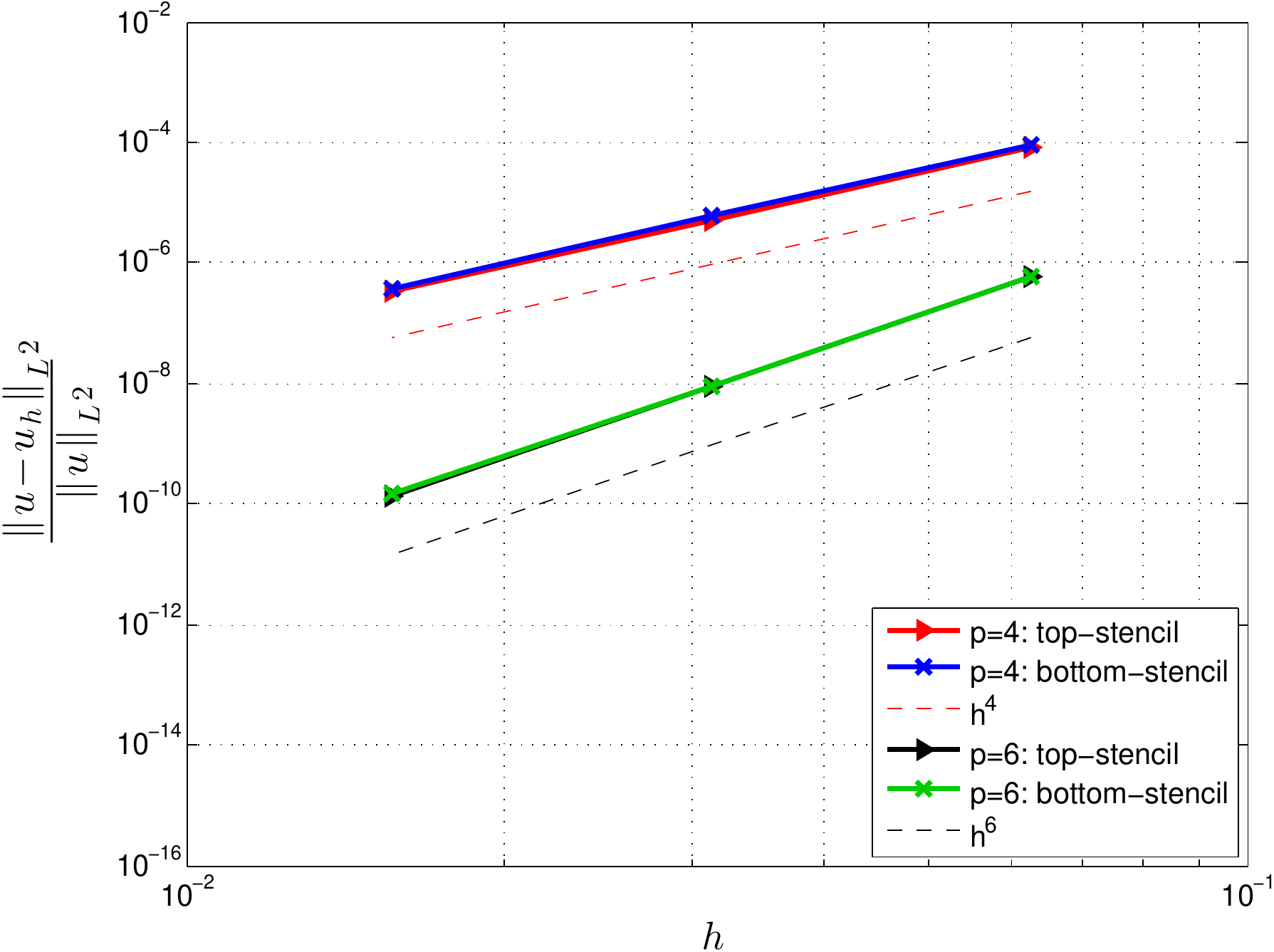}
\includegraphics[width=0.44\linewidth]{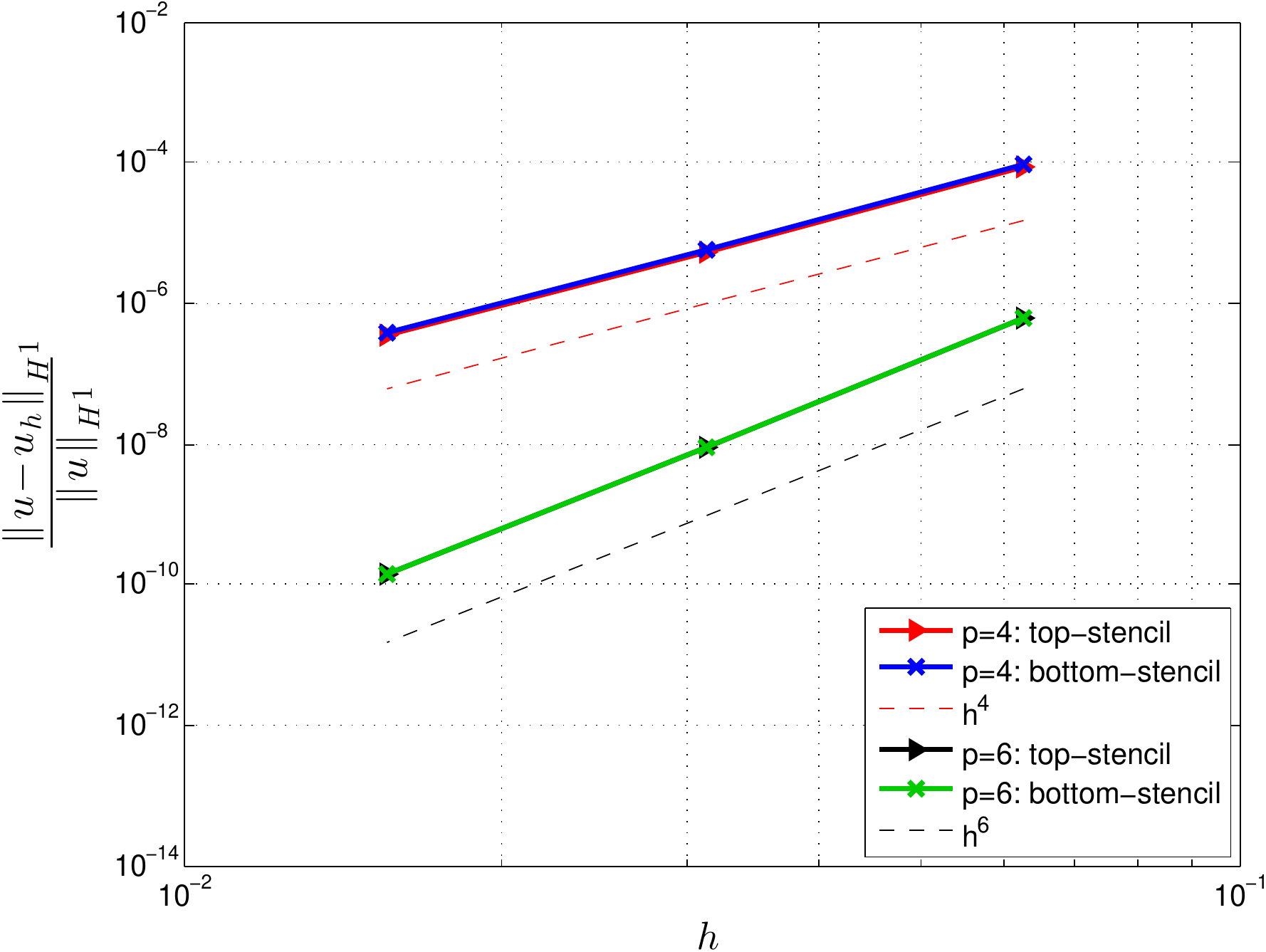}
\caption{$L^{2}$ and $H^{1}$ error plot for  C-CSP periodic problem (even
  $p$),  with the stencils depicted in Figure
  \ref{periodic_ccsp_even}.}
\label{CSPperiodicerror-even}
\end{center}
\end{figure}

A natural question arises: why can't we achieve optimal convergence
when even degrees B-splines are considered?  The answer is not yet
clear. As we explained in the previous sections, 
the rationale behind C-CSP, as well as LS-SP and C-ASP, is to
try to obtain the same solution delivered by the Galerkin method by
imposing the residual to be zero at the superconvergent
points, which are supposedly close to the true zeros of the Galerkin
residual.
 However, as we discussed in Section \ref{sec:sppoints}, we do not have access to the
precise location of the superconvergent points, and instead we use
``surrogate'' superconvergent points that do not approximate well
the zeros of the Galerkin residual everywhere in the domain.  We do not see
however any qualitative difference between the odd and even case
other than in the central element 
({although we did not perform a quantitative analysis of this issue}). 
Furthermore, it is not clear why the
C-ASP points would have poorer approximation properties than the C-CSP ones: 
in other words, the points that would be selected by the C-ASP seem as good as those
that would be selected by the C-CSP as for what concerns being close
to the zeros of the Galerkin residual.  It could be that the local
symmetry of the C-CSP points distribution leads to some  error
cancellation in the collocation system.

We continue by testing the C-CSP method on the Dirichlet Problem \ref{prob:dirichlet1d}
with $a_0=a_1=0$ and $f(x) = \pi ^2\sin(\pi x)$, whose exact solution
is $u(x)=\sin(\pi x)$, and show the corresponding results in Figure $\ref{CSP1d}$.
As in the previous case, the order of convergence is $p+1$ in $L^{2}$ norm and $p$ in $H^{1}$ norm. 
In order to compare the four methods we presented (C-CSP, ASP,
LS-SP and Greville collocation) against the Galerkin solver, 
we show in Figure $\ref{Comparison1d}$ a comparison 
of the convergence of the $L^2$-error obtained when solving the
Dirichlet problem above with B-splines of degree $p=3$.
The plot highlights that C-CSP, although converging with optimal order,
shows an error one order of magnitude larger than Galerkin, 
while LS-SP converges essentially to the same solution of the Galerkin
method. It should be observed however that the computational cost of LS-SP
is significantly higher than C-CSP, not only because of
the number of points where the residual needs to be evaluated 
(about $2^d$ times more that C-CSP in $d$ dimensions) 
but also for the the higher condition number of the resulting system of
linear equations. 

\begin{figure}[tp]
\begin{center}
\includegraphics[width=0.44\textwidth]{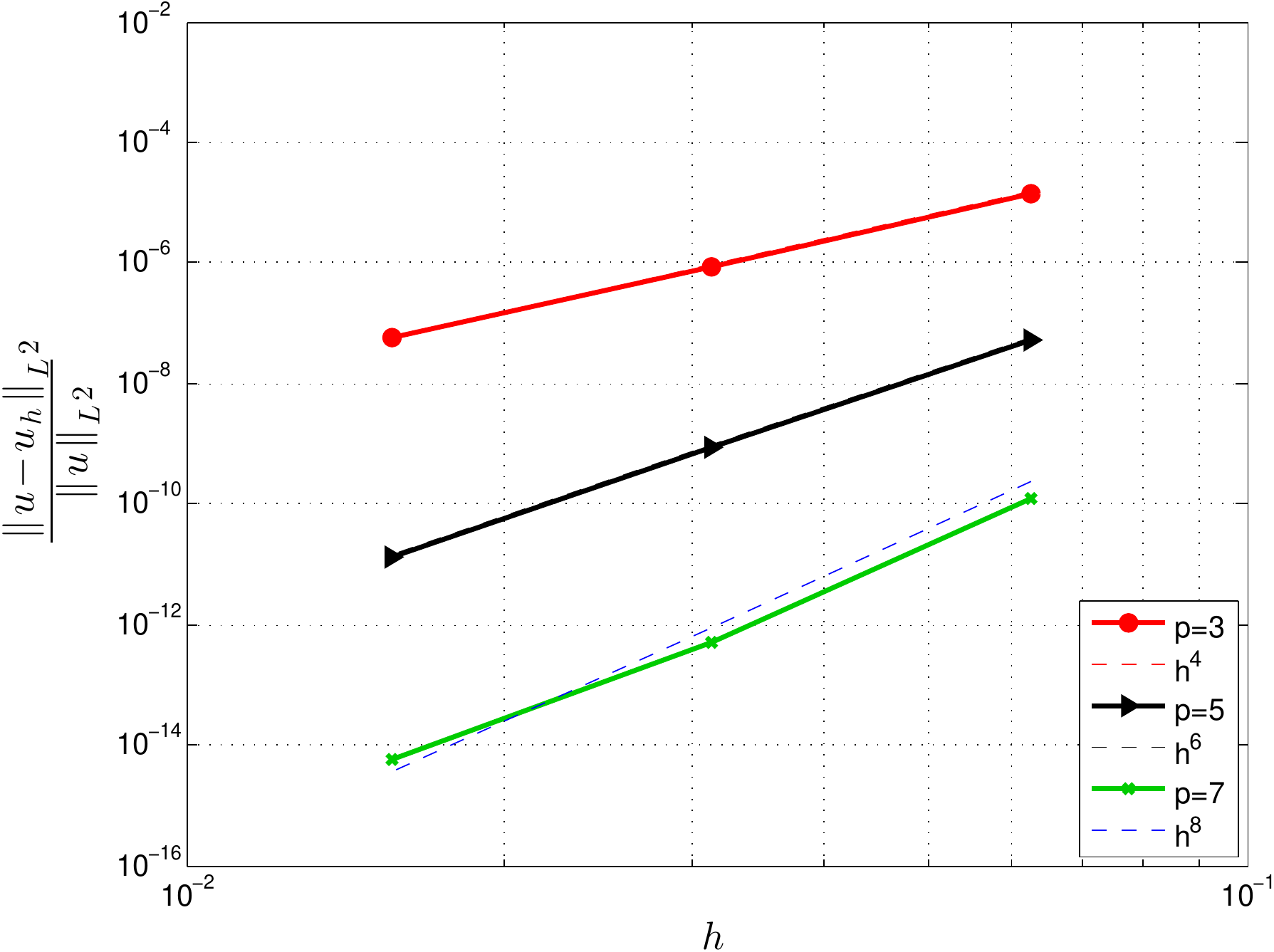}
\includegraphics[width=0.44\textwidth]{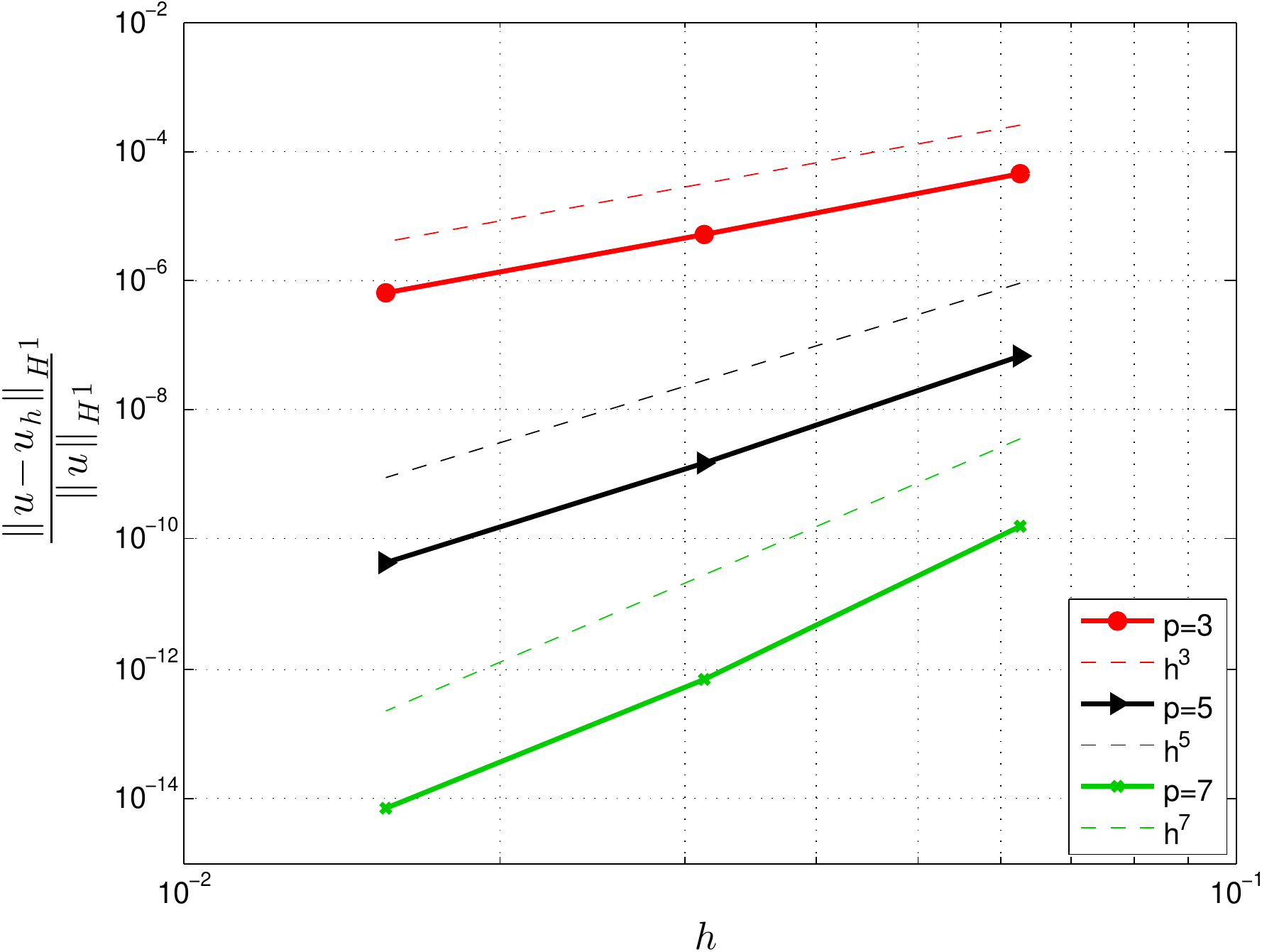}
\caption{$L^{2}$ and $H^{1}$  error plot: C-CSP Dirichlet problem.}
\label{CSP1d}
\end{center}
\end{figure}

\begin{figure}[tp]
\begin{center}
\includegraphics[width=0.44\textwidth]{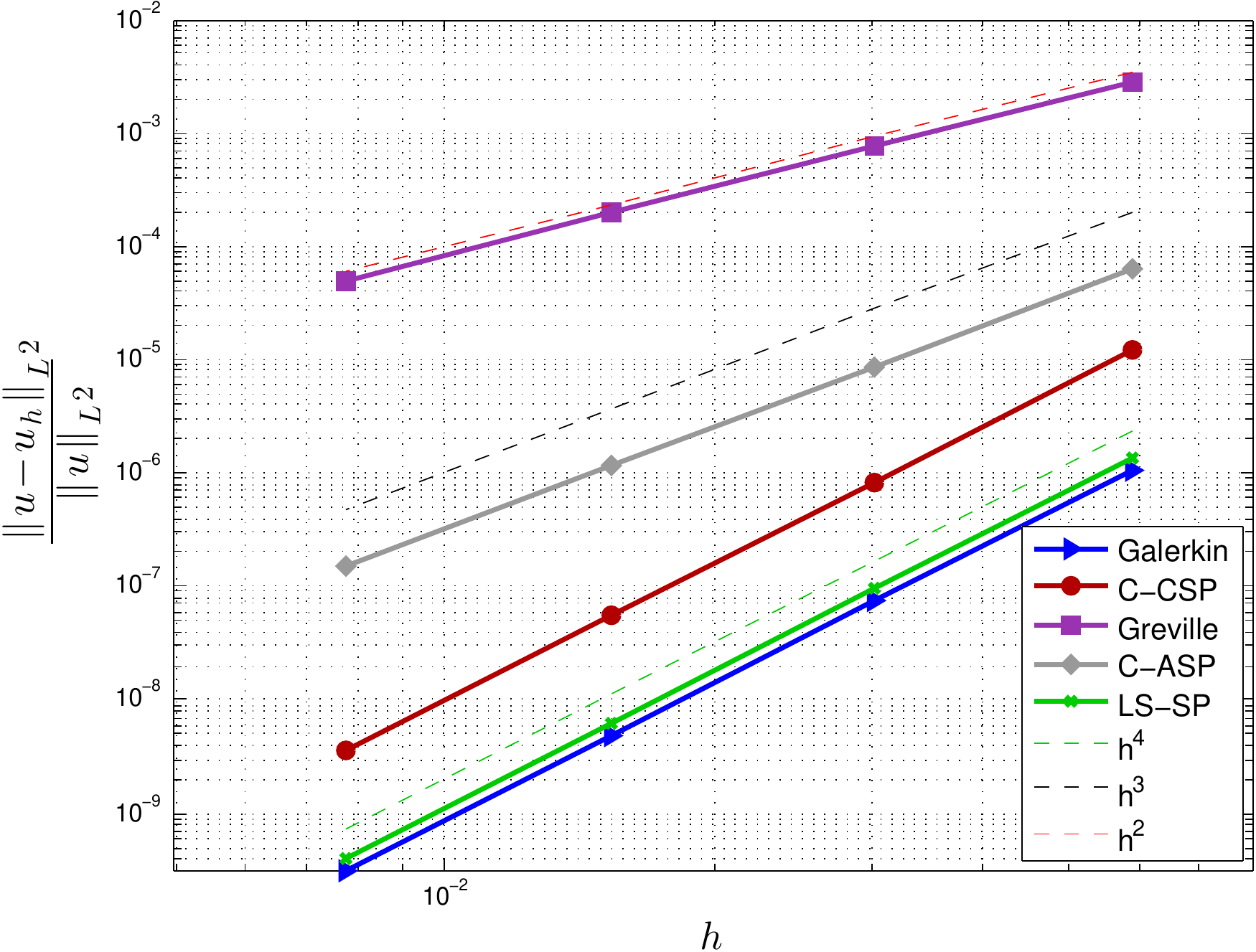}
\caption{Comparison of convergence of $L^2$ error norms for the Dirichlet problem
for different methods.}
\label{Comparison1d}
\end{center}
\end{figure}

We also investigate the robustness of the method with respect to
perturbations of the knot vector. To this end, 
we perturb the internal knots of each equispaced open knot vector considered in the 
convergence analysis by} randomly chosen quantities, i.e. 
we replace each internal knot $\xi_i$ by $\tilde{\xi}_i = \xi_i + \frac{1}{10 n_{el}} X_i$, 
where $X_i$ are independent random numbers $X_i \in [-1,1]$. {We remark that the random quantities $X_i$
are generated at every refinement step for each node of the knot vector.
 The scaling factor $\frac{1}{10 n_{el}}$ prevents 
knot clashes and furthemore the resulting knot
vectors are  quasi-uniform, but the local symmetry of the
mesh is lost for all elements of the mesh. This is expected to have an
influence on the location of the superconvergent points, but we nonetheless select
the collocation points following the element-wise construction for the uniform mesh case.  
The error plots are shown in Figures \ref{fig:rand.knots.L2} and \ref{fig:rand.knots.H1}. 
We note that we loose the optimal rates of convergence we observed in the previous tests: 
the order of convergence is $p$ for both the $L^2$ and $H^1$ error norm, i.e., optimal for the $H^1$ error norm 
and only one-order suboptimal for the $L^2$ one.

We also verify the influence of the differential operator and of the boundary conditions, by 
considering non-null $a_0$ and $a_1$ in the Dirichlet Problem \ref{prob:dirichlet1d}, 
as well as Neumann-Neumann and Neumann-Dirichlet boundary conditions.
We performed several tests and the results obtained were identical; therefore, 
we report here only one representative example, {see \cite{montardini:MSthesis} for additional numerical results}. 
In detail, we consider 
$a_1(x)=x$, $a_0=1$  and $f(x) = x(e^x\sin(\pi x) + \pi e^x\cos(\pi x)) - 2\pi e^x\cos(\pi x) + \pi^2 e^x\sin(\pi x)$,
whose exact solution is $u(x)=\sin(\pi x)e^x$. 
The results of the test we performed are shown in 
Figures \ref{fig:adv.reac:L2} and \ref{fig:adv.reac:H1}: 
the order of convergence is still optimal:  $p+1$ 
for the $L^2$ error norm, and $p$ for the $H^1$ norm.
We can then conclude that the C-CSP method seems to be robust 
with respect the form of the elliptic operator.

\begin{figure}[tp]
 \centering
 \subfloat[][random-perturbation test, $L^2$ error.\label{fig:rand.knots.L2}]
   {\includegraphics[width=.44\textwidth]{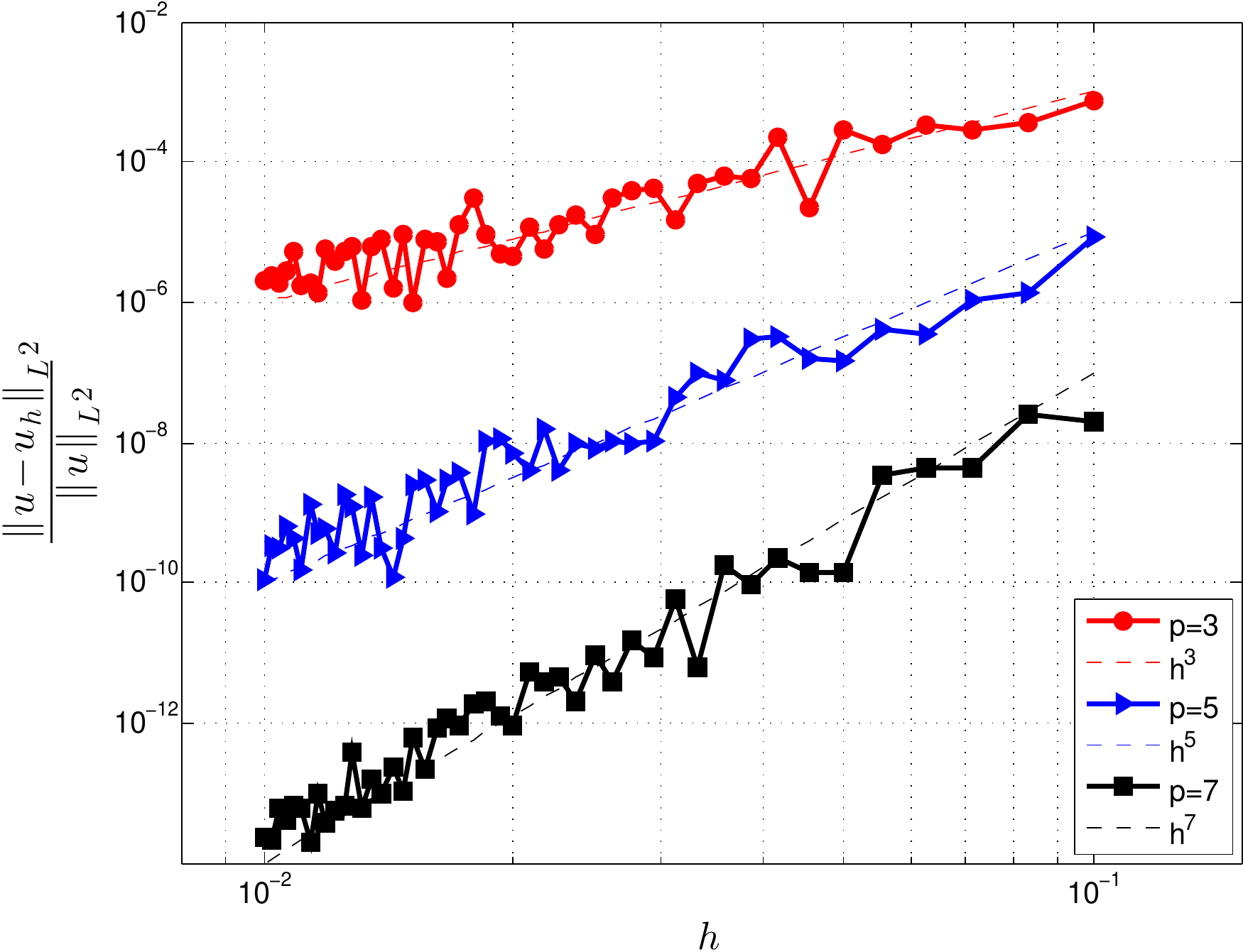}}\quad
 \subfloat[][random-perturbation test, $H^1$ error.\label{fig:rand.knots.H1}]
   {\includegraphics[width=.44\textwidth]{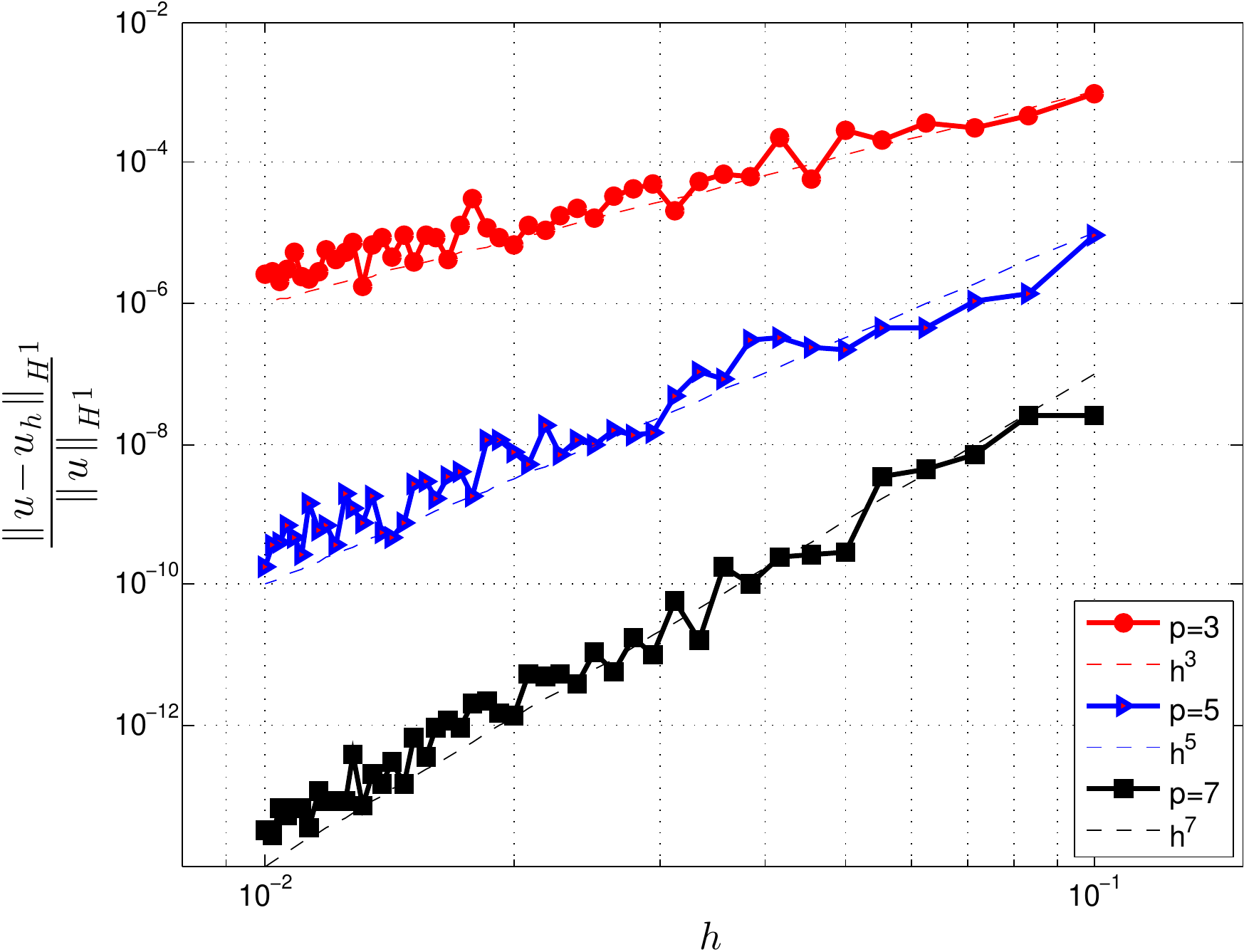}}\quad
 \subfloat[][advection-reaction test, $L^2$ error.\label{fig:adv.reac:L2}]
   {\includegraphics[width=.44\textwidth]{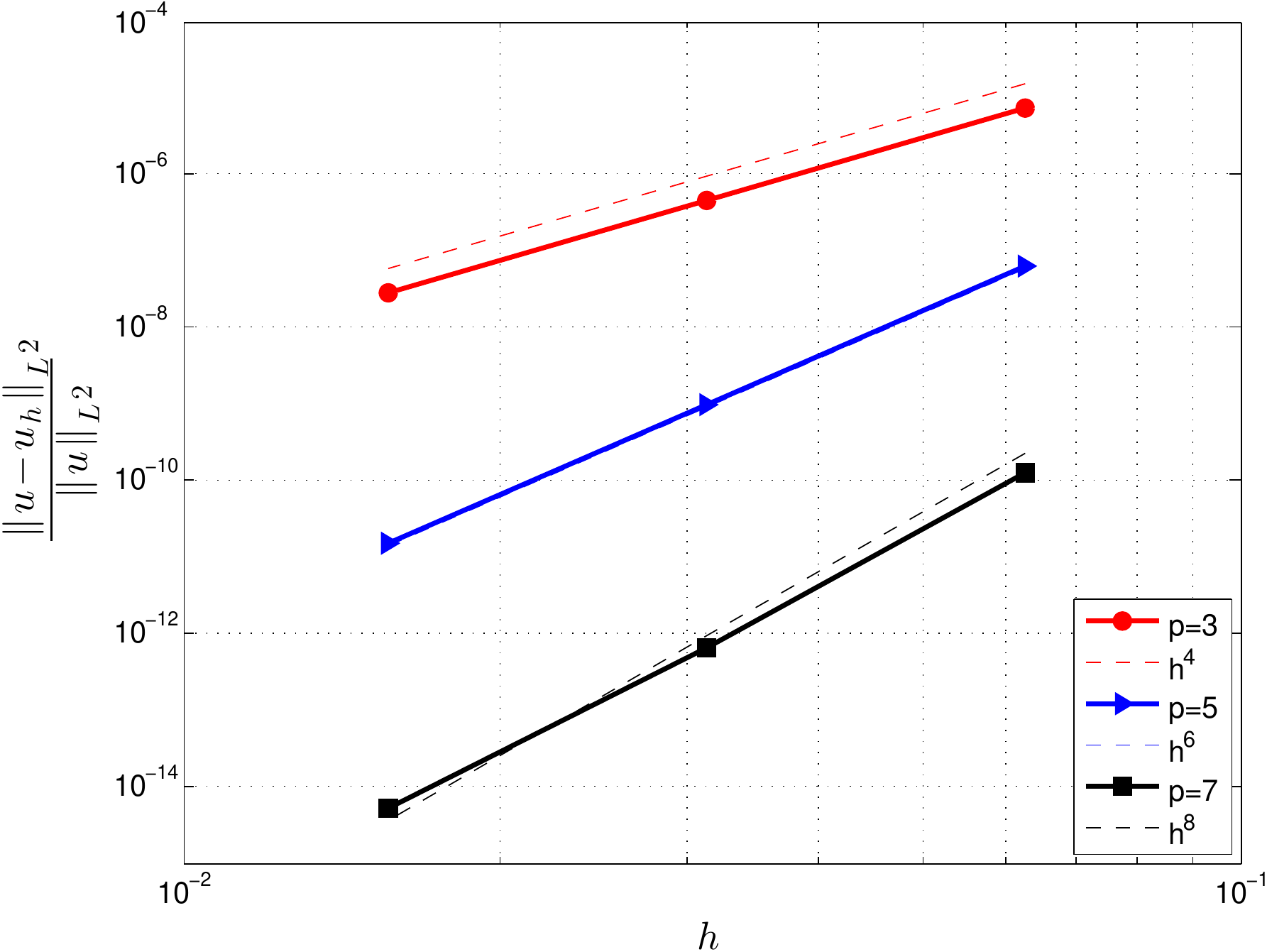}}\quad
 \subfloat[][advection-reaction test, $H^1$ error.\label{fig:adv.reac:H1}]
   {\includegraphics[width=.44\textwidth]{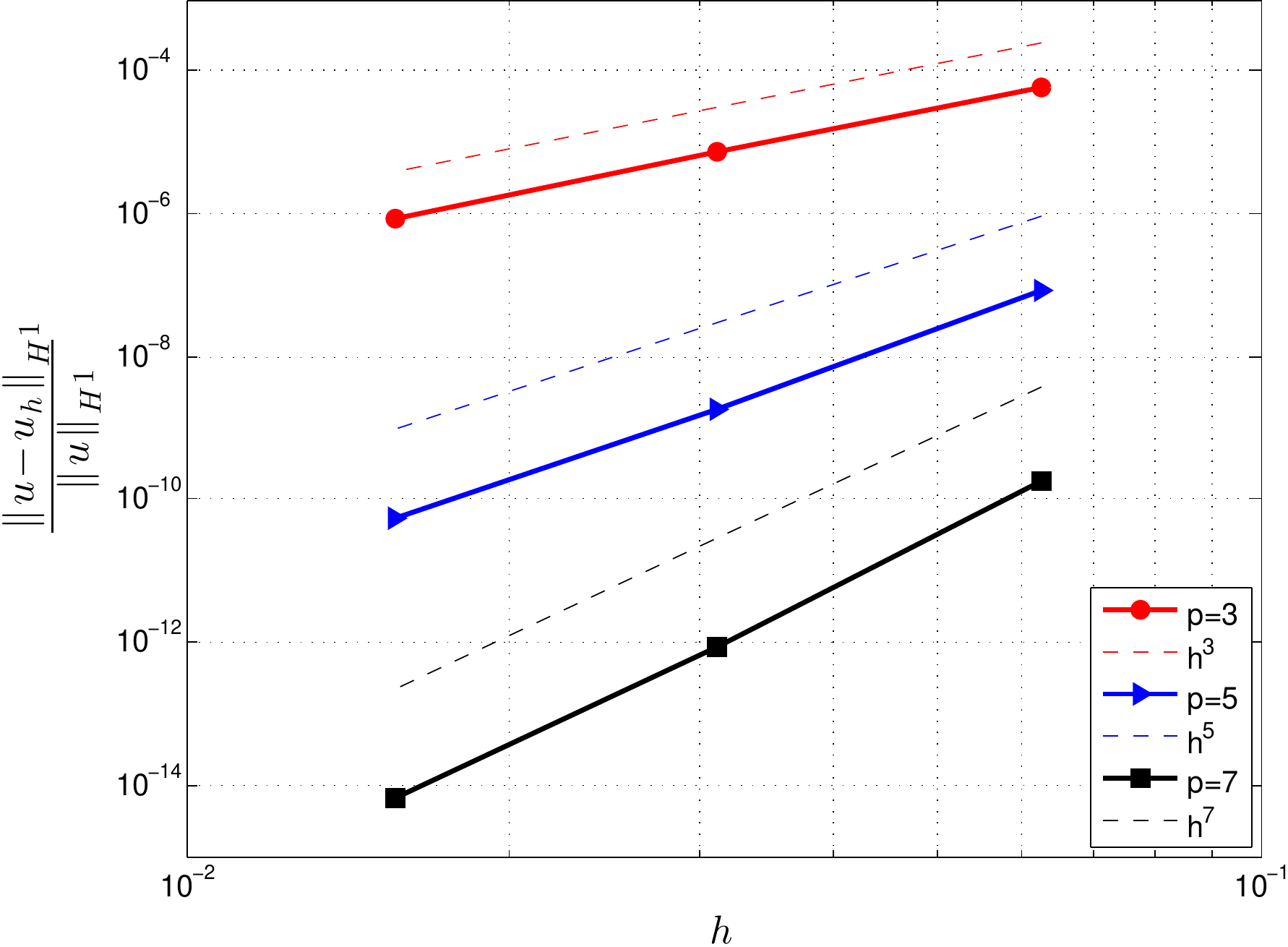}}
 \caption{Robustness test for C-CSP with respect to 
   perturbation of the knot vectors (plots \ref{fig:rand.knots.L2} and \ref{fig:rand.knots.H1}) 
   and changes of the differential operator (plots \ref{fig:adv.reac:L2} and \ref{fig:adv.reac:H1}). }
 \label{Tests}
 \end{figure}

Finally, we present two examples of two-dimensional Dirichlet Problem \ref{prob:dirichlet} solved by C-CSP.
In the first one we consider as computational domain,
$\Omega$, the quarter of annulus in Figure \ref{ring}
(for which NURBS functions have to be employed), and $f(x,y)$ is chosen
such that the exact solution is $u(x,y)=-(x^2+y^2-1)(x^2+y^2-4)xy^2$.
In Figures \ref{2Dringcspall:p3}, \ref{2Dringcspall:p5} and \ref{2Dringcspall:p7} we show the convergence
plots of the $L^2$ and $H^1$ errors for the C-CSP and Galerkin methods for odd degree NURBS $p=3,5,7$. 
The observed orders of convergence are as expected: optimal order in $L^2$ norm and $H^1$ norm.
We remark that for $p=3$ the $H^1$ norms of the corresponding errors
obtained by the Galerkin and C-CSP methods are very close, 
as can be seen in Figure \ref{2Dringcspall:p3}.

 \begin{figure}[tp]
   \hspace{0.5cm}
   \qquad \subfloat[][Ring domain. \label{ring}]  {\includegraphics[height=0.35\textwidth]{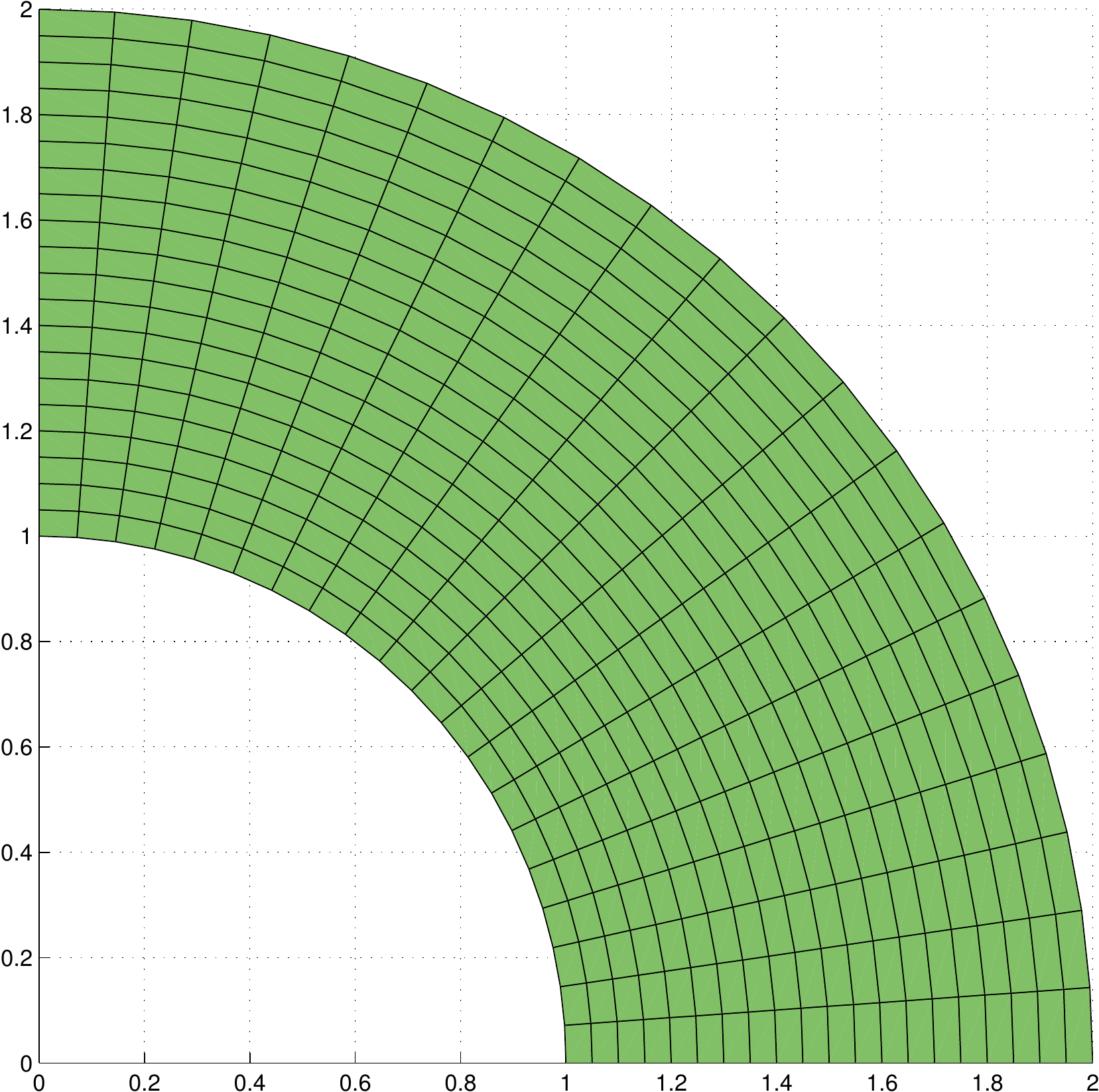}}\qquad \qquad 
   \subfloat[][Case $p=3$. \label{2Dringcspall:p3}]{\includegraphics[width=.49\textwidth]{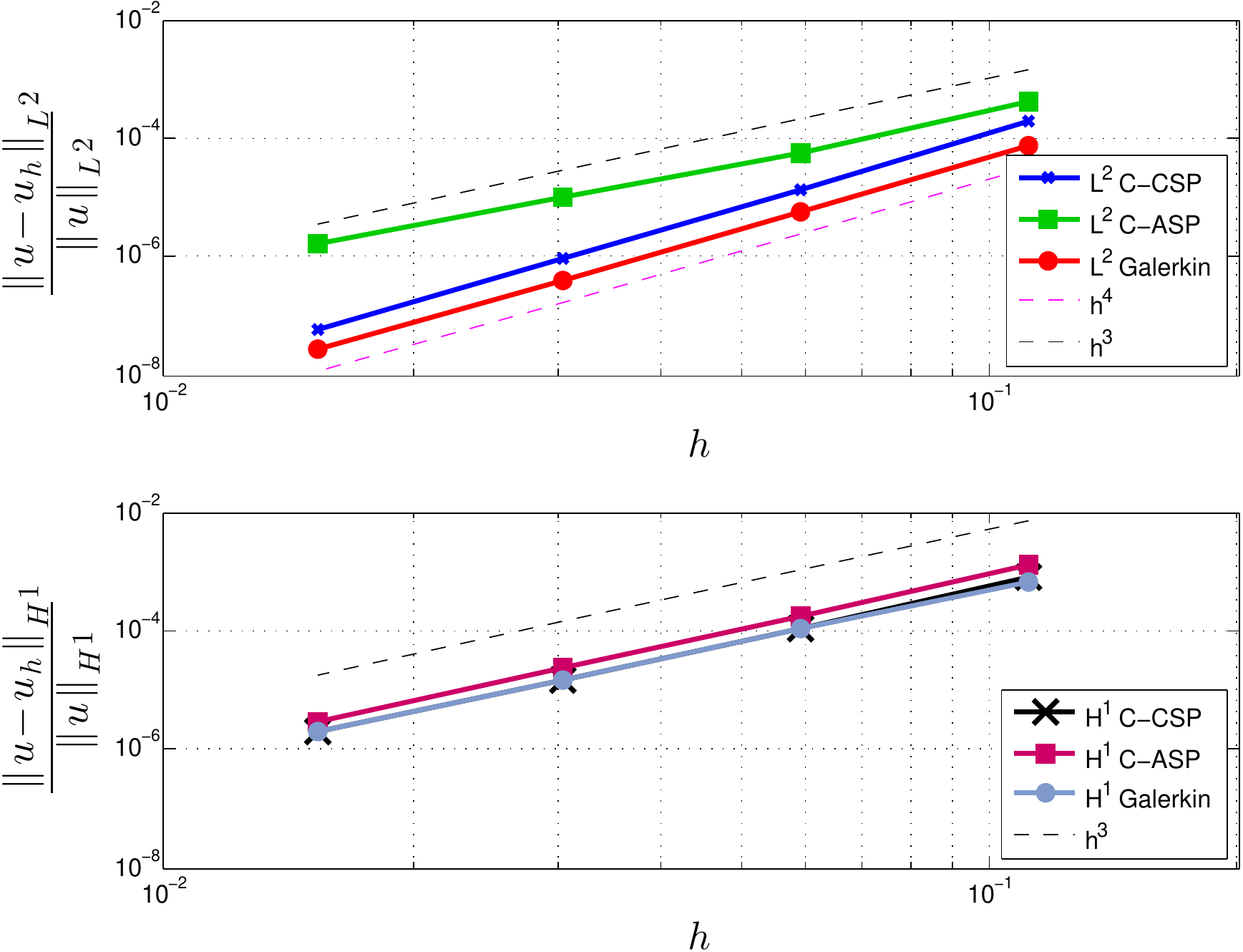}} \\
   \subfloat[][Case $p=5$. \label{2Dringcspall:p5}]{\includegraphics[width=.49\textwidth]{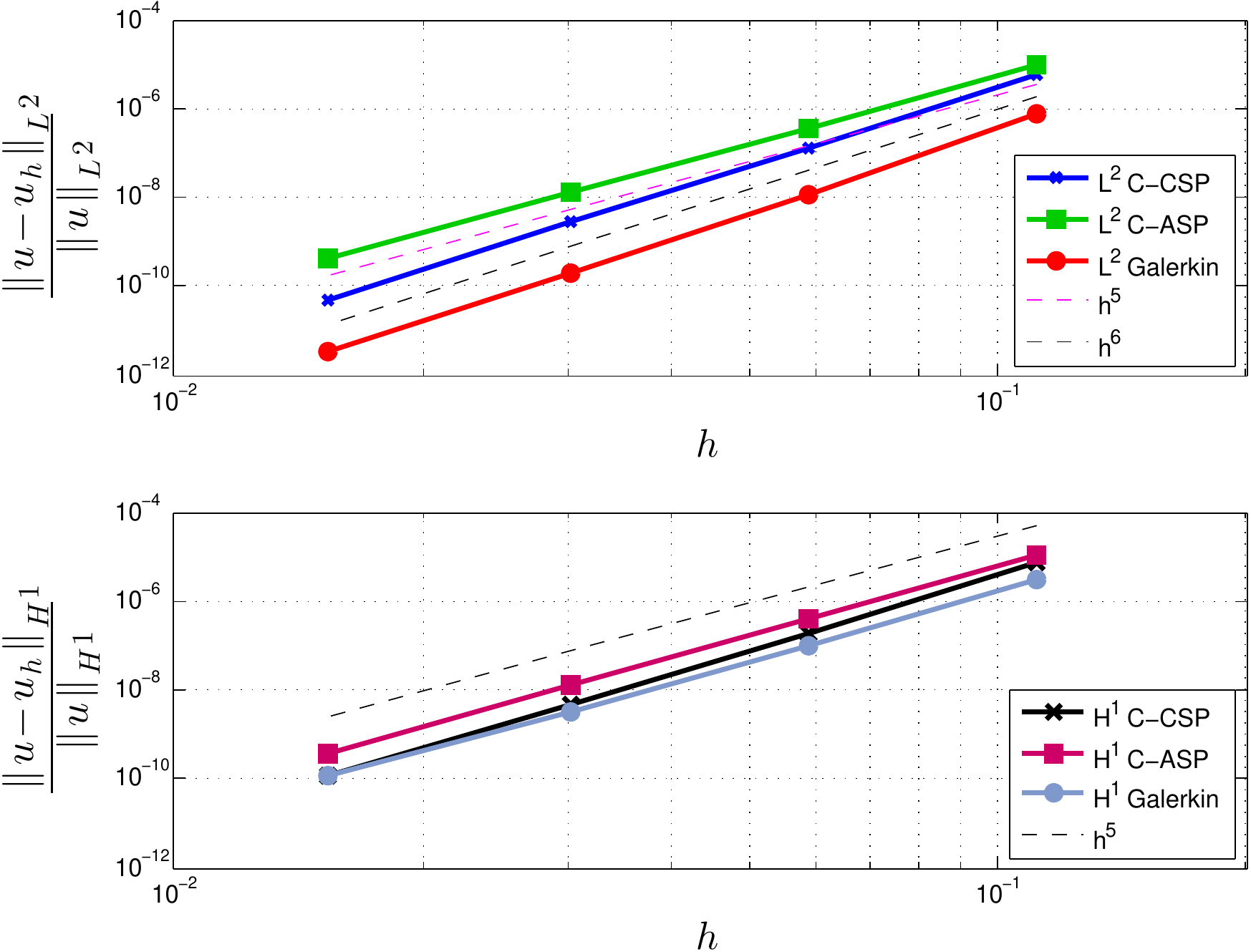}}\quad
   \subfloat[][Case $p=7$. \label{2Dringcspall:p7}]{\includegraphics[width=.49\textwidth]{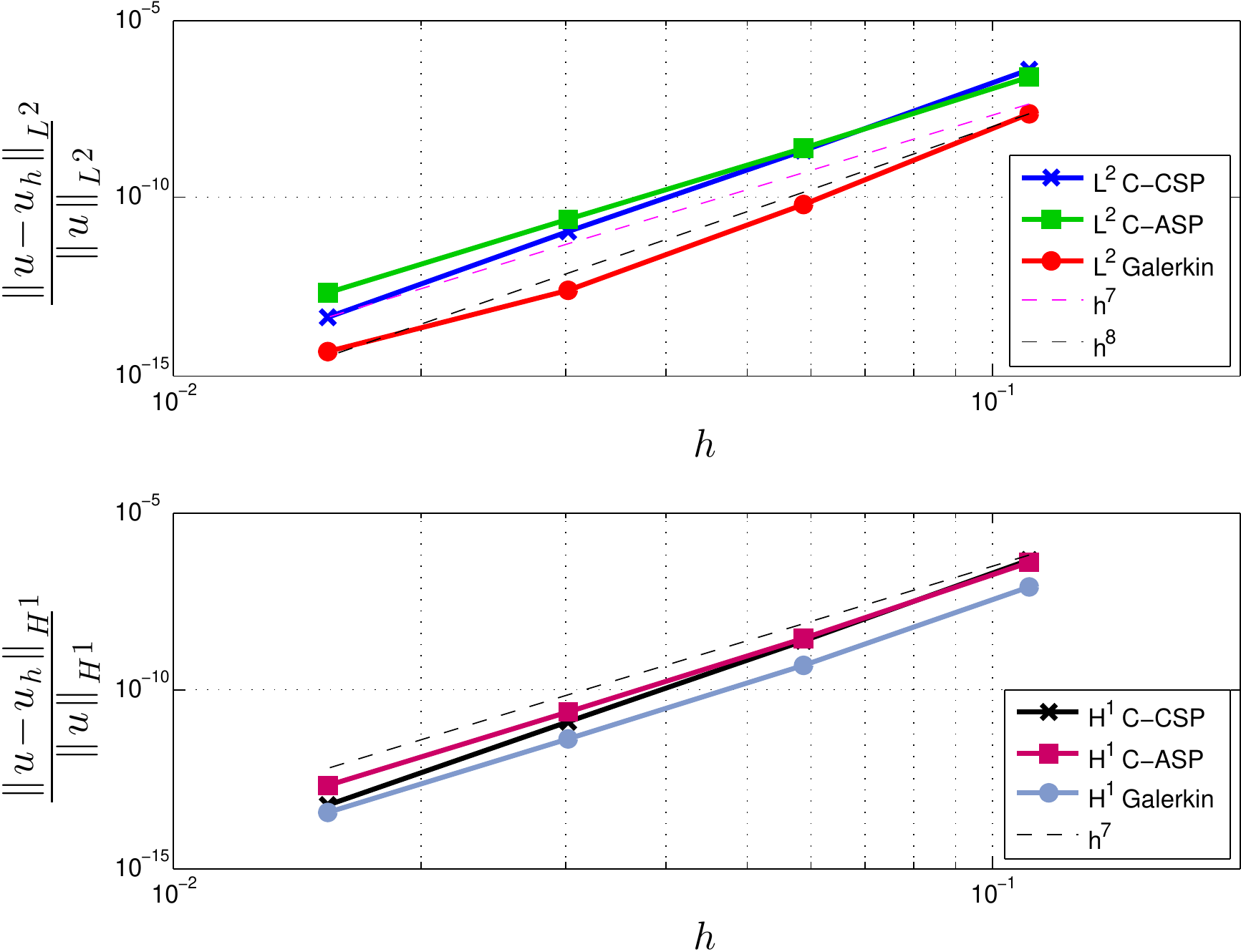}}
   \caption{Ring domain and $L^2$ and $H^1$ convergence of approximations with $p=3,5,7$. 
     }
   \label{2Dringcspall}
 \end{figure}

In the second two-dimensional example, we let $\Omega$
be the rhombus with vertices $(0,0)$, $(\frac{1}{4}, 1)$,
$(1, \frac{1}{4})$ and $(\frac{5}{4}, \frac{5}{4})$, represented
in Figure \ref{rhombus} (its parametrization is bilinear but not
orthogonal, as in the previous example), and $f$ is such that the exact solution is
$u(x,y)=\sin\left(\frac{4}{15}\pi(y-4x)\right)\sin\left(\frac{16}{15}\pi\left(\frac{x}{4}-y\right)\right)(x^3+y^3)$.
The corresponding errors are shown in Figures \ref{2Drhombusall:p3}, \ref{2Drhombusall:p5} 
and \ref{2Drhombusall:p7}, and the same observations as before hold.
Note however that the gap between the C-CSP and the Galerkin solution is
larger than in the previous example, especially for the $L^2$ error.
Moreover, for $p=5$ and $p=7$ the convergence is still in its preasymptotic regime. 

 \begin{figure}[tp]
   \hspace{0.5cm}
  \qquad   \subfloat[][Rhombus domain \label{rhombus}]{\includegraphics[height=0.35\textwidth]{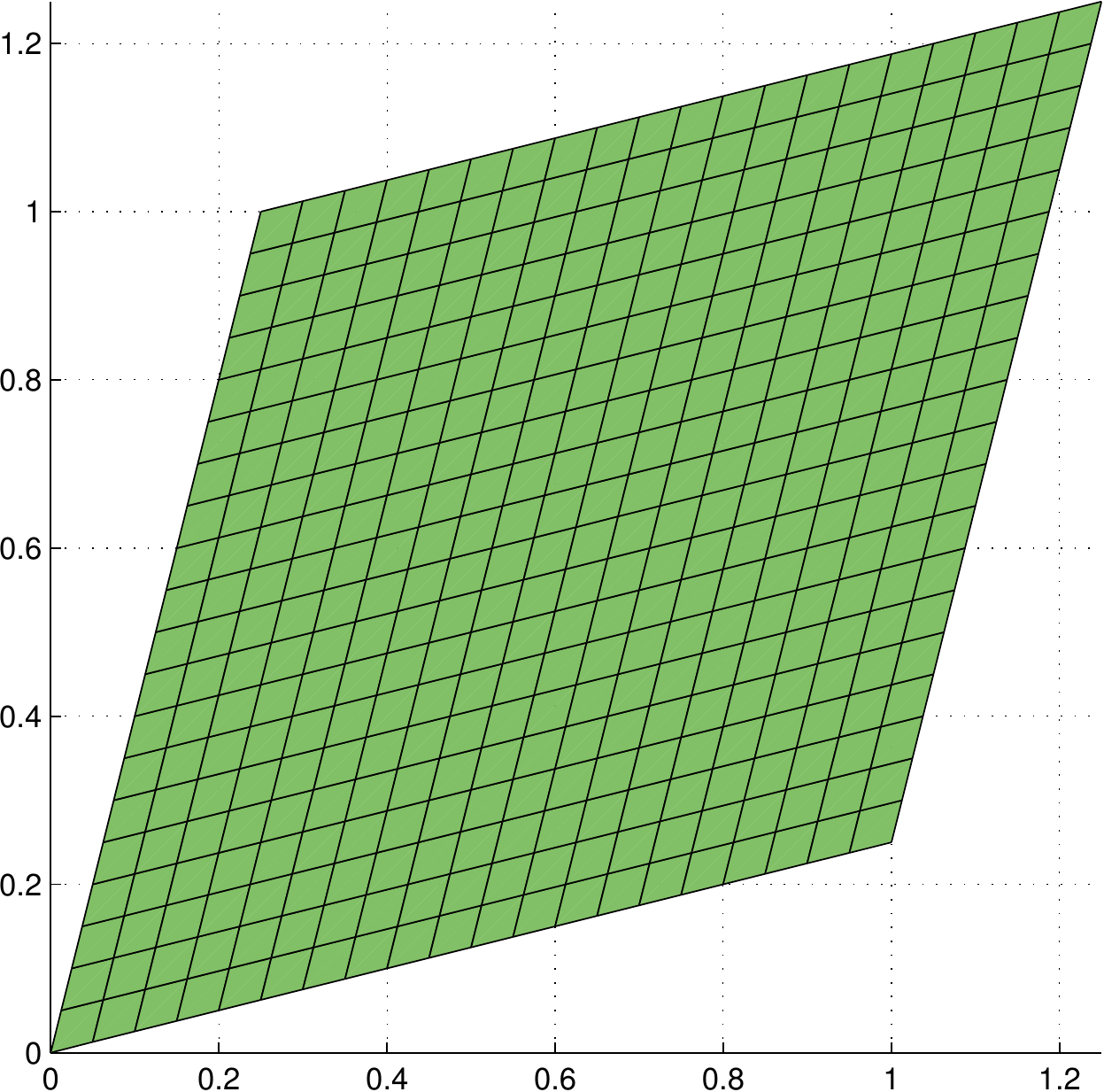}} \qquad \qquad 
   \subfloat[][Case $p=3$. \label{2Drhombusall:p3}]{\includegraphics[width=.49\textwidth]{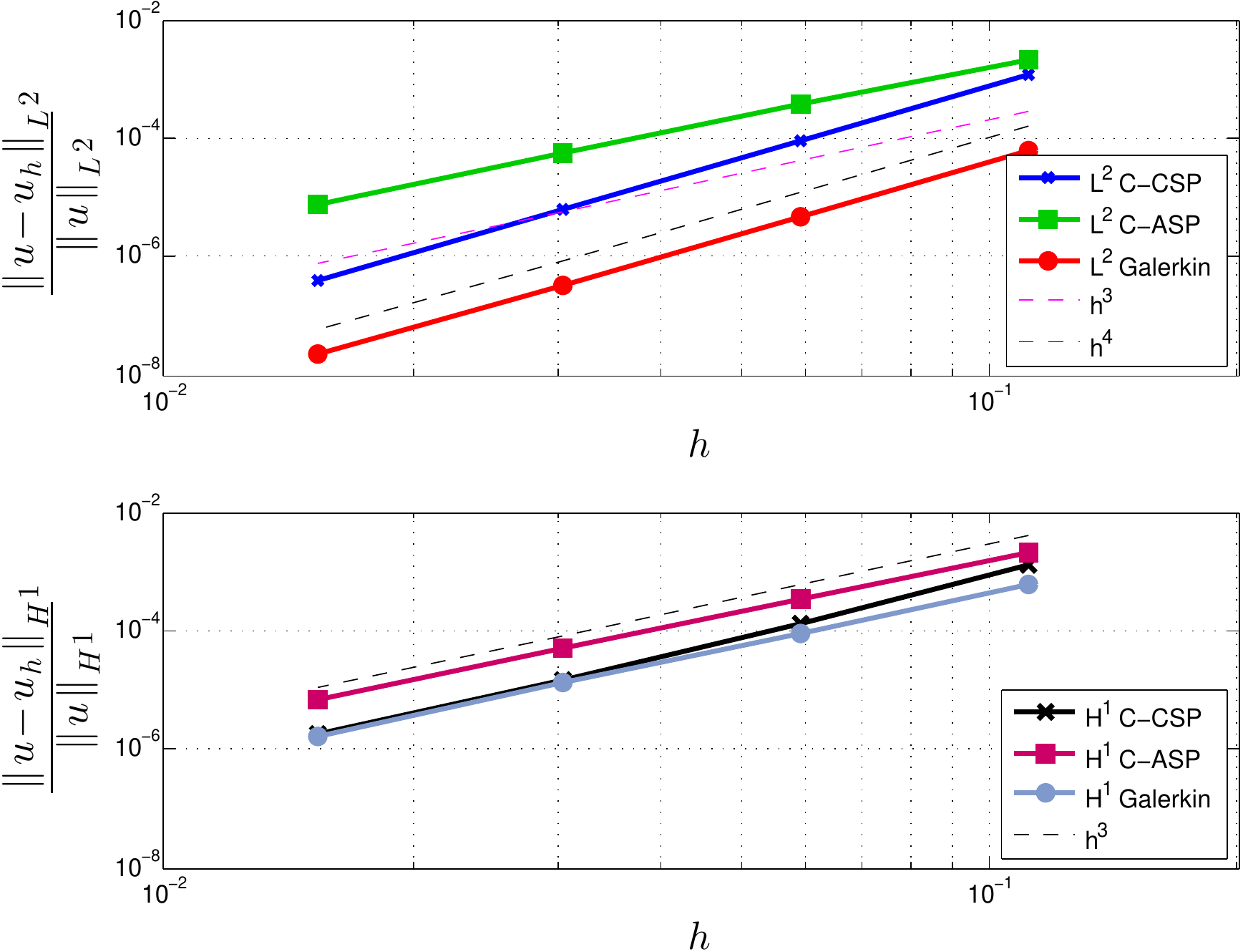}}\\
   \subfloat[][Case $p=5$. \label{2Drhombusall:p5}]{\includegraphics[width=.49\textwidth]{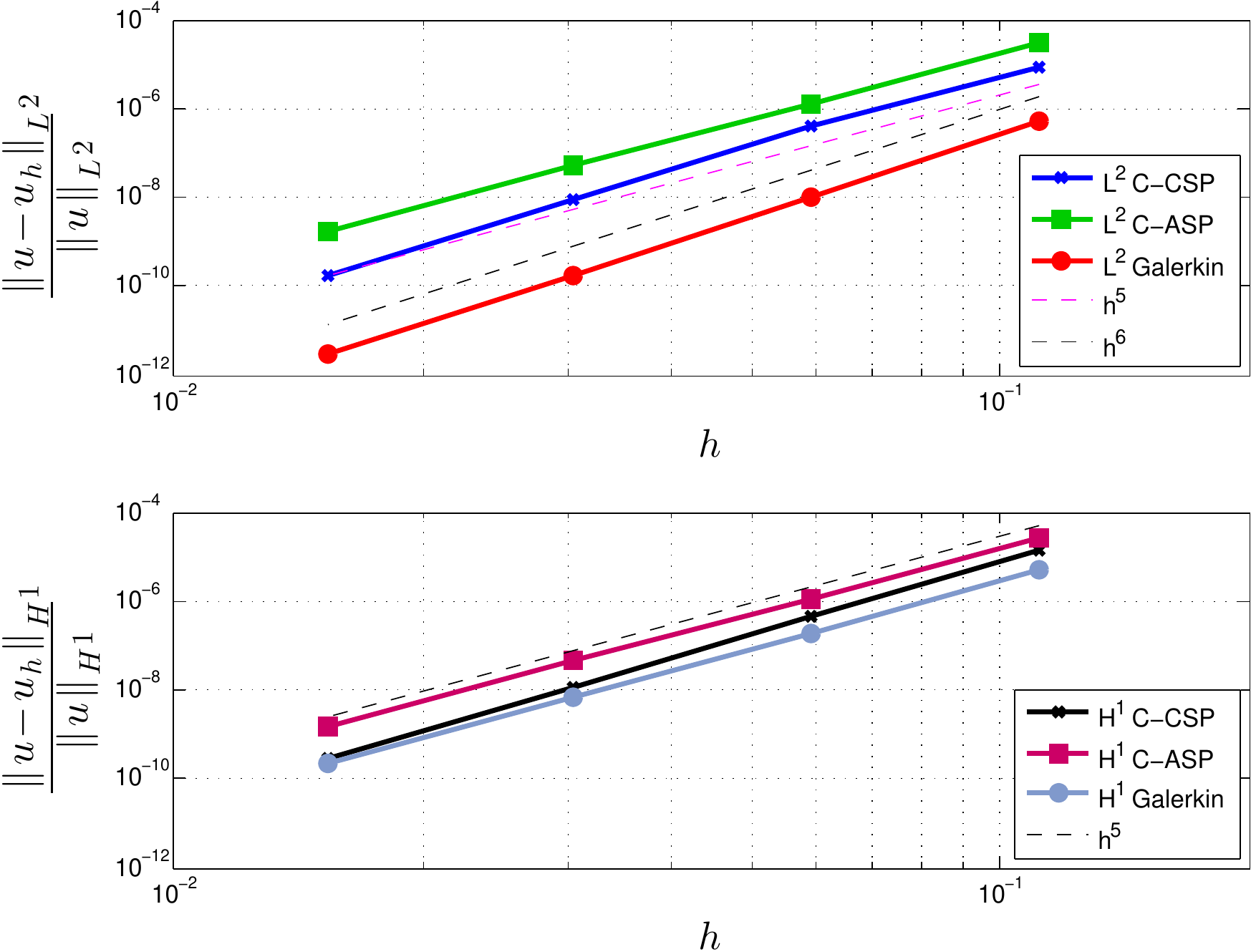}}
   \subfloat[][Case $p=7$. \label{2Drhombusall:p7}]{\includegraphics[width=.49\textwidth]{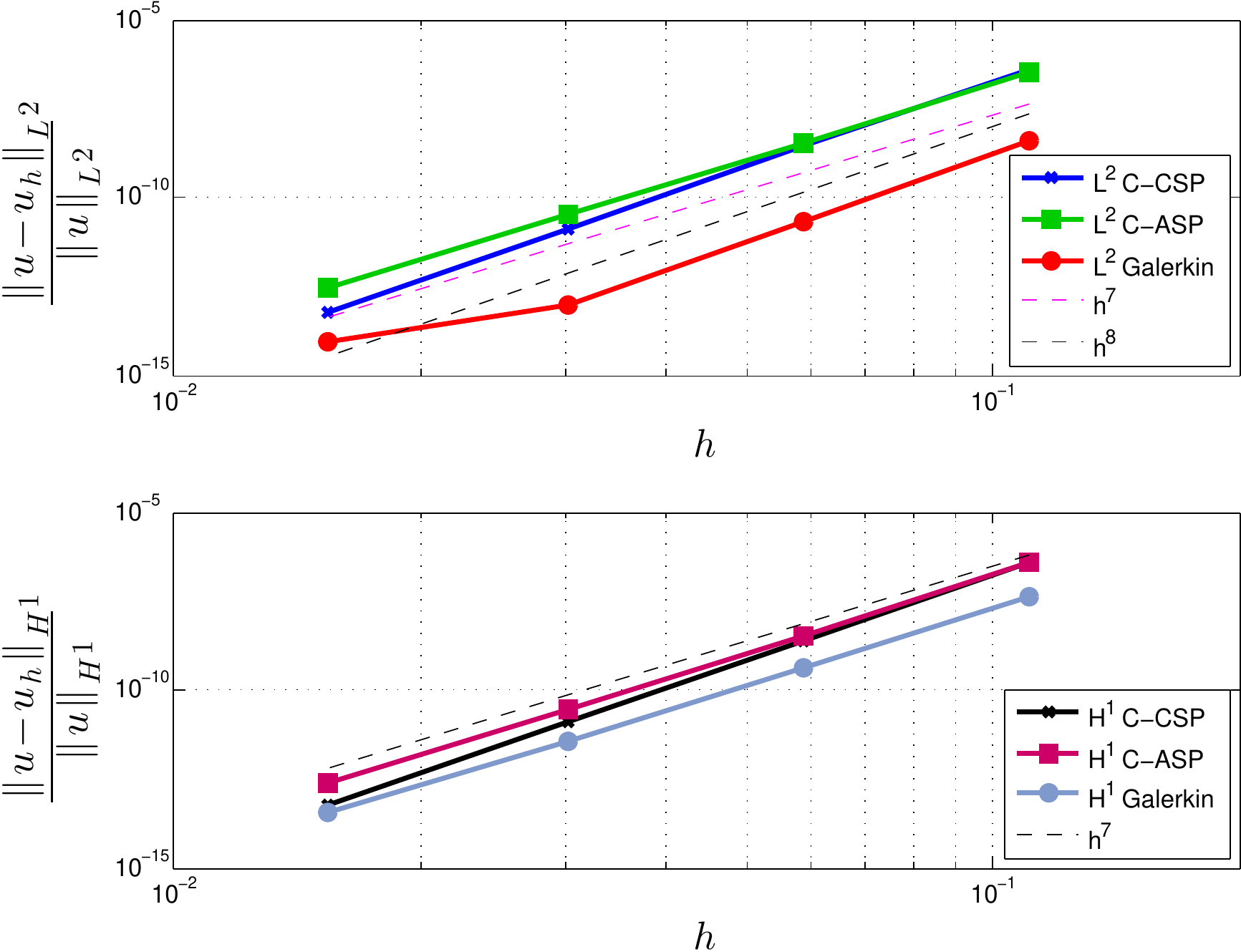}}
   \caption{Rhombus domain and $L^2$ and $H^1$ convergence of approximations with $p=3,5,7$. 
     }
   \label{2Drhombusall}
 \end{figure}

\section{Conclusion}\label{sec:conclusion}
In this paper we have proposed an isogeometric collocation method
based on the superconvergent Galerkin points, in the framework of
\cite{gomez2016variational}.   
Our guiding criterion  for the subset of the superconvergent points is
however different, and consists in picking up 
clusters of points, in such a way to 
{obtain a collocation scheme symmetric at element-scale}.  
This choice allows to recover optimal convergence rates for odd-degrees splines/NURBS,
 without  ``oversampling'' the domain as in least-square approach proposed by
\cite{anitescu2015isogeometric}.
Moreover, the order of convergence of the
$L^\infty$ norm of the error is the same of $L^2$-norm in all experiments we
have performed (not shown in the paper for the sake of brevity).

The preliminary numerical campaign on one- and two-di\-men\-sional tests
that we performed suggests that
the method is robust with respect to isogeometric mapping of the
domain, while  perturbations of the knot vector may reduce the accuracy of the method.  
A rigorous mathematical explanation for the convergence
behavior observed for the proposed method, and for the other collocation
methods based on the Galerkin superconvergent points, is not
available yet and will be the target of our future efforts.

\section*{Acknowledgements}
The authors were partially supported by the European Research Council
through the FP7 ERC Consolidator Grant n.616563 \emph{HIGEOM},
by European Union's Horizon 2020 research and innovation program 
through the grant no. 680448 CAxMan 
and by the Italian MIUR through the PRIN  
``Metodologie innovative nella modellistica differenziale numerica''.  
This support is gratefully acknowledged.

\def\cprime{$'$}

\end{document}